\documentclass{article}
\usepackage{amssymb}
\usepackage{amsmath}
\usepackage{latexsym}
\usepackage{pictex}
\begin{document}

\bibliographystyle{plain}
\baselineskip=15pt
\newtheorem{lemma}{Lemma}[section]
\newtheorem{theorem}[lemma]{Theorem}
\newtheorem{prop}[lemma]{Proposition}
\newtheorem{cor}[lemma]{Corollary}
\newtheorem{definition}[lemma]{Definition}
\newtheorem{definitions}[lemma]{Definitions}
\newtheorem{remark}[lemma]{Remark}
\newtheorem{conj}[lemma]{Conjecture}
\renewcommand{\arraystretch}{2}
\newcommand{\TE}{\widetilde{E}}
\newcommand{\TF}{\widetilde{F}}
\newcommand{\TP}{\widetilde{P}}
\newcommand{\C}{{\mathbb{C}}}
\newcommand{\N}{{\mathbb{N}}}
\newcommand{\Q}{{\mathbb{Q}}}
\newcommand{\R}{{\mathbb{R}}}
\newcommand{\Z}{{\mathbb{Z}}}
\newcommand{\B}{{\mathbf{B}}}
\newcommand{\CA}{{\mathcal{A}}}
\newcommand{\CL}{{\mathcal{L}}}
\newcommand{\0}{-\,}
\newcommand{\dimv}{{\underline{\dim}}}
\newcommand{\cdi}{\mbox{CD$(\mathbf{i})$}}
\newcommand{\cdii}{\mbox{CD$(\mathbf{i}')$}}
\newcommand{\Hom}{{{\rm Hom}_{kQ}}}
\newcommand{\Ext}{{{\rm Ext}^1_{kQ}}}
\begin{center}
\baselineskip=15pt
{\LARGE Canonical basis linearity regions arising from quiver representations}

\vspace*{1cm}

{\Large Robert Marsh} \\
\vspace*{0.5cm}
{\em Department of Mathematics and Computer Science, University of Leicester,
University Road, Leicester LE1 7RH, England} \\
\vspace*{0.5cm}
{E-mail: R.Marsh@mcs.le.ac.uk}

\vspace*{0.3cm}
{\Large Markus Reineke} \\
\vspace*{0.5cm}
{\em BUGH Wuppertal,
         Gau\ss stra\ss e 20,
         D-42097 Wuppertal,
         Germany} \\
\vspace*{0.5cm}
{E-mail: reineke@math.uni-wuppertal.de}

\baselineskip=15pt

{\bf Abstract} \\
\parbox[t]{5.15in}{
In this paper we show that there is a link between the combinatorics of the
canonical basis of a quantized enveloping algebra and the monomial bases
of the second author~\cite{reineke3}
arising from representations of quivers. We prove that
some reparametrization functions of the canonical basis, arising from
the link between Lusztig's approach to the canonical basis and the string
parametrization of the canonical basis, are given on a large region
of linearity by linear functions arising from these monomial
basis for a quantized enveloping algebra.

} \\ \ \\
\end{center}

\section{Introduction}

Let $U=U_q({\mathbf{g}})$ be the quantum group associated to a semisimple
Lie algebra $\mathbf{g}$ of rank $n$. The negative part $U^-$ of $U$ has
a canonical basis $\mathbf{B}$ with favourable properties
(see Kashiwara~\cite{kash2} and Lusztig~\cite[\S14.4.6]{lusztig6}).
For example, via action on highest weight vectors it gives rise to bases
for all the finite-dimensional irreducible highest weight $U$-modules.

The reparametrization functions of the canonical basis studied in this paper
arise from considering two parametrizations of the canonical basis, both
dependent on the choice of a reduced expression for the longest
element $w_0$ in the corresponding Weyl group compatible with a quiver.
The first parametrization
arises from correspondence between the canonical basis and a
basis of PBW type, and the second is the string parametrization
(see~\cite[\S2]{bz1}, the end of Section $2$ in~\cite{nz1} and~\cite{kash4}).
The reparametrization functions are useful in relating these two
parametrizations, and in~\cite{bz3} they were shown to have applications to
the description of the tensor product multiplicities of simple
$\mathbf{g}$-modules in an approach involving totally positive varieties.
They were also shown to be closely connected with the description
of the canonical basis in~\cite{cartermarsh} and~\cite{marsh9}.

Our approach to these functions differs from that of Berenstein
and Zelevinsky~\cite{bz3} in that we make the linearity of these functions
on certain regions explicit, and that we show how they are linked with the
representation theory of quivers. In particular, we show that, 
for reduced expressions compatible with quivers, they can be
described, on certain cones (which we call degeneration cones), as
linear functions arising from the monomial basis of $U^-$ of
the second author~\cite{reineke3} which was found using properties of quiver representations,
and in particular, degenerations of representations.

We show that the degeneration cones contain the Lusztig cones
of~\cite{lusztig7}. Thus the reparametrization functions mentioned above,
when restricted to Lusztig cones, can be described using quiver theory. Lusztig cones have been
shown to have many interpretations and connections. For example, they are
used to describe regular functions on a reduced real double Bruhat cell of
the corresponding algebraic group~\cite{zelevinsky}, they have links with
primitive elements in the dual canonical basis (this can be seen
using~\cite{bz1}) and therefore with the representation theory
of affine Hecke algebras~\cite{lnt}, they are known to correspond to
regions of linearity of the Lusztig reparametrization functions and tight
monomials (see~\cite{cartermarsh}), and they can be described
using the homological algebra of representations of quivers~\cite{bedard2}.

The structure of the paper is as follows. In \S2 we recall the
parametrizations of the canonical basis we will need, and in \S3 we recall
the connection between reduced expressions for the longest word in the Weyl
group and quivers. In \S4 we recall the monomial basis and corresponding
linear functions defined in~\cite{reineke3}, as well as showing
that the functions we will use are invertible. In \S5 we recall the Lusztig
cones of~\cite{lusztig7}, and in \S6 we define the degeneration cones
referred to above. In \S7 we define a set of cones in the PBW parametrization,
which we show in \S9 correspond the Lusztig cones in appropriate string
parametrizations; in \S8 we show that each cone in \S7 is contained in the
corresponding degeneration cone. Finally, in section \S10 we show that
the linear functions arising in~\cite{reineke3} coincide with the
reparametrization functions on the degeneration cones, using the description
of the operators of Kashiwara given in~\cite{reineke1}.

\section{Parametrizations of the canonical basis}

Let $\mathbf{g}$ be the simple Lie algebra over $\C$ of type $A_n$ and $U$ be
the
quantized enveloping algebra of $\mathbf{g}$. Then $U$ is a $\Q(v)$-algebra
generated by the elements $E_i$, $F_i$, $K_{\mu}$, $i\in \{1,2,\ldots ,n\}$,
$\mu\in Q$, the root
lattice of $\mathbf{g}$. Let $U^+$ be the subalgebra generated by the $E_i$ and
$U^-$ the subalgebra generated by the $F_i$.

Let $W$ be the Weyl group of $\mathbf{g}$. It has a unique element $w_0$
of maximal length. For each reduced expression $\mathbf{i}$ for $w_0$ there
are two
parametrizations of the canonical basis $\B$ for $U^-$. The first arises from
Lusztig's approach to the canonical basis~\cite[\S14.4.6]{lusztig6}, and
the second arises from Kashiwara's approach~\cite{kash2}.

\noindent {\bf Lusztig's Approach} \\
There is an $\Q$-algebra automorphism of $U$ which takes each $E_i$ to
$F_i$, $F_i$ to $E_i$, $K_{\mu}$ to $K_{-\mu}$ and $v$ to $v^{-1}$.
We use this automorphism to transfer Lusztig's
definition of the canonical basis in~\cite[\S3]{lusztig2} to $U^-$.

Let $T_i$, $i=1,2,\ldots ,n$, be the automorphism of $U$ as
in~\cite[\S1.3]{lusztig1} given by:

$$T_i(E_j)=\left\{ \begin{array}{cc}
-F_jK_j, & {\rm if\ }i=j, \\
E_j,     & {\rm if\ }|i-j|>1 \\
-E_iE_j+v^{-1}E_jE_i & {\rm if\ }|i-j|=1
\end{array} \right.$$

$$T_i(F_j)=\left\{ \begin{array}{cc}
-K_j^{-1}E_j, & {\rm if\ }i=j, \\
F_j,     & {\rm if\ }|i-j|>1 \\
-F_jF_i+vF_iF_j & {\rm if\ }|i-j|=1
\end{array} \right.$$

$$T_i(K_{\mu})=K_{\mu-\langle \mu,\alpha_i\rangle h_i},{\rm \ for\ }\mu\in
Q,$$ where the $\alpha_i$ are the simple roots and the $h_i$ are the
simple coroots of $\mathbf{g}$.

For each $i$, let $r_i$ be the automorphism of $U$ which fixes $E_j$ and
$F_j$ for $j=i$ or $|i-j|>1$ and fixes $K_{\mu}$ for all $\mu$,
and which takes $E_j$ to $-E_j$ and $F_j$ to $-F_j$ if $|i-j|=1$.
Let $T''_{i,-1}=T_ir_i$ be the automorphism of $U$ as
in~\cite[\S37.1.3]{lusztig6}.
Let ${\mathbf{c}}\in {\N}^N$, where $N=\ell (w_0)$, and
$\mathbf{i}$ be a reduced expression for $w_0$. Let
$$F_{\mathbf{i}}^{\mathbf{c}}:=F_{i_1}^{(c_1)}T''_{i_1,-1}(F_{i_2}^{(c_2)})\cdots
T''_{i_1,-1}T''_{i_2,-1}\cdots T''_{i_{N-1},-1}(F_{i_N}^{(c_N)}).$$
Define $B_{\mathbf{i}}=\{F_{\mathbf{i}}^{\mathbf{c}}\,:\, {\mathbf{c}}\in
{\N}^N\}$.
Then $B_{\mathbf{i}}$ is
the basis of PBW-type corresponding to the reduced expression $\mathbf{i}$.
Note that, if the reduced expression $\mathbf{i}$ is adapted to a
quiver in the sense of~\cite{lusztig2}, then this basis can also be
constructed using the Hall algebra approach of~\cite{ringel}.
Let\,\, $\mathbf{\bar{\ }}\,$ be the ${\Q}$-algebra automorphism from $U$ to
$U$ taking $E_i$ to $E_i$, $F_i$ to $F_i$, and $K_{\mu}$ to $K_{-\mu}$,
for each $i\in [1,n]$ and $\mu\in Q$, and $v$ to $v^{-1}$.
Lusztig proves the following result in~\cite[\S\S2.3, 3.2]{lusztig2}.

\begin{theorem} (Lusztig) \\
The $Z[v]$-span $\CL$ of $B_{\mathbf{i}}$ is independent of $\mathbf{i}$.
Let $\pi:{{\CL}}\rightarrow {{\CL}}/{v{\CL}}$ be the natural
projection. The image $\pi(B_{\mathbf{i}})$ is also independent
of $\mathbf{i}$; we denote it by $B$. The restriction of $\pi$ to
${{\CL}}\cap \overline{{\CL}}$ is an isomorphism of $\mathbb{Z}$-modules
$\pi_1:{{\CL}}\cap \overline{{\CL}}\rightarrow {{\CL}}/{v{\CL}}$.
Also $\B=\pi_1^{-1}(B)$ is a
${\Q}(v)$-basis of $U^-$, which is the canonical basis of $U^-$.
\end{theorem}

Lusztig's theorem provides us with a parametrization of ${\B}$, dependent on
$\mathbf{i}$. If $b\in {\B}$, we write $\phi_{\mathbf{i}}(b)={\mathbf{c}}$,
where ${\mathbf{c}}\in {\N}^N$ satisfies $b\equiv F_{\mathbf{i}}^{\mathbf{c}}
\mod v{{\CL}}$.
Note that $\phi_{\mathbf{i}}$ is a bijection.

For reduced expressions $\mathbf{j}$ and $\mathbf{j'}$ for $w_0$,
Lusztig defines in~\cite[\S2.6]{lusztig2} a reparametrization
function $R_{\mathbf{j}}^{\mathbf{j}'}=\phi_{\mathbf{j}'}
\phi_{\mathbf{j}}^{-1}\,:\,{\N}^N \rightarrow {\N}^N$. This function
was shown by Lusztig to be piecewise linear and its regions of linearity were
shown to have significance for the canonical basis, in the sense that
elements $b$ of the canonical basis with $\phi_{\mathbf{j}}(b)$ in the same
region of linearity of $R_{\mathbf{j}}^{{\mathbf{j}}'}$ often have
similar form. For example, this can be seen from the explicit descriptions of 
the canonical basis of type $A_2$ and $A_3$, as
computed by Lusztig in~\cite{lusztig2} and by Xi in~\cite{xi}, respectively. 
More evidence for the importance of these regions is their connection~\cite{bz1} with the multiplicativity properties of dual canonical bases.

\noindent {\bf The string parametrization}

Let $\widetilde{E}_i$ and $\widetilde{F}_i$ be the Kashiwara operators on
$U^-$ as defined in~\cite[\S3.5]{kash2}. Let ${\CA}\subseteq {\Q}(v)$ be
the subring of elements regular at $v=0$, and let ${\CL}'$ be the
${\CA}$-lattice spanned by arbitrary products
$\widetilde{F}_{j_1}\widetilde{F}_{j_2}\cdots \widetilde{F}_{j_m}\cdot 1$ in
$U^-$.
We denote the set of all such elements by $S$. The following results
were proved by Kashiwara in~\cite{kash2}.

\begin{theorem} \label{kashiwara} (Kashiwara) \\
(i) Let $\pi':{{\CL}'}\rightarrow {{\CL}'}/{v{{\CL}'}}$ be the natural
projection, and let $B'=\pi'(S)$. Then $B'$ is a ${\Q}$-basis of
${{\CL}'}/{v{{\CL}'}}$ (the crystal basis). \\
(ii) Furthermore,
$\widetilde{E}_i$ and $\widetilde{F}_i$ each preserve ${\CL}'$ and thus act on
${{\CL}'}/{v{{\CL}'}}$. They satisfy $\widetilde{E}_i(B')\subseteq
B'\cup\{0\}$ and
$\widetilde{F}_i(B')\subseteq B'$. Also for $b,b'\in B'$ we have
$\widetilde{F}_ib=b'$,
if and only if $\widetilde{E}_ib'=b$. \\
(iii) For each $b\in B'$, there is a unique element $\widetilde{b}\in {{\CL}'}
\cap \overline{{\CL}'}$ such that $\pi'(\widetilde{b})=b$. The set of elements
$\{\widetilde{b}\,:\, b\in B'\}$ forms a basis of $U^-$, the {\em global
crystal basis} of $U^-$.
\end{theorem}

It was shown by Lusztig~\cite[2.3]{lusztig3} that the global crystal basis
of Kashiwara coincides with the canonical basis of $U^-$.

There is a parametrization of ${\B}$ arising from Kashiwara's approach,
again dependent on a reduced expression $\mathbf{i}$ for $w_0$.
Let $\mathbf{i}=(i_1,i_2,\ldots ,i_k)$ and $b\in {\B}$.
Let $a_1$ be maximal such that $\widetilde{E}_{i_1}^{a_1}b\not\equiv 0 \mod
v{{\CL}'}$;
let $a_2$ be maximal such that
$\widetilde{E}_{i_2}^{a_2}\widetilde{E}_{i_1}^{a_1}b\not\equiv 0 \mod
v{{\CL}'}$,
and so on, so that
$a_N$ is maximal such that
$$\widetilde{E}_{i_N}^{a_N}\widetilde{E}_{i_{N-1}}^{a_{N-1}}
\cdots \widetilde{E}_{i_2}^{a_2}\widetilde{E}_{i_1}^{a_1}b\not\equiv 0 \mod
v{{\CL}'}.$$
Let $\mathbf{a}=(a_1,a_2,\ldots ,a_N)$. We write
$\psi_{\mathbf{i}}(b)=\mathbf{a}$.
This is the crystal string of $b$ --- see~\cite[\S2]{bz1}
and the end of Section $2$ in~\cite{nz1}; see also~\cite{kash4}.
It is known that $\psi_{\mathbf{i}}(b)$ uniquely determines $b\in {\B}$
(see~\cite[\S2.5]{nz1}).
We have $b\equiv \TF_{i_1}^{a_1}\TF_{i_2}^{a_2}\cdots \TF_{i_N}^{a_N}\cdot 1
\mod v{\CL}'$.
The image of $\psi_{\mathbf{i}}$ is a cone which appears
in~\cite{bz1}.
We shall call this the {\em string cone} $X_{st}(\mathbf{i})=
\psi_{\mathbf{i}}(\B)$.

We next consider a function which links the string parametrization with
the parametrization arising from Lusztig's approach.
For reduced expressions $\mathbf{i}$ and $\mathbf{j}$ for
$w_0$, consider the maps
$$X_{st}(\mathbf{i})\xrightarrow[\psi_{\mathbf{i}}^{-1}]{}\B
\xrightarrow[\phi_{\mathbf{j}}]{} \N^N$$
We define $S_{\mathbf{i}}^{\mathbf{j}}=\phi_{\mathbf{j}}
\psi_{\mathbf{i}}^{-1}\,:\,X_{st}(\mathbf{i})\rightarrow {\N}^N$,
a reparametrization function. This function has appeared in,
for example,~\cite{bz3}.

{\bf Remark on Notation:}

We shall use the following notation convention. A cone $C$ which
is to be regarded as a subset of $X_{st}$, i.e. to be thought of as a set
of strings, will be given the subscript ``st'' to denote this, thus $C_{st}$.
A cone $C$ which is to be regarded as a subset of $\mathbb{N}^N$, and is
to be regarded as a set of PBW parameters for the canonical basis, will be
denoted with the subscript ``PBW'', thus $C_{PBW}$. In each case, it will be
clear from the context which reduced expression for $w_0$ is being used.

\section{Reduced Expressions Compatible with Quivers} \label{compatible}

We now fix the type of $\mathbf{g}$ once and for all to be type $A_n$.
However, a lot of the results we use hold in a greater generality, and we
also believe that the results proved here also hold in greater generality
(at least for the simply-laced Dynkin case).
In this section we shall recall an explicit description of reduced expressions
compatible with quivers in type $A_n$, which we shall need in studying such
expressions. We shall write a
reduced expression $s_{i_1}s_{i_2}\cdots s_{i_N}$ for $w_0$ as the $N$-tuple
$\mathbf{i}=(i_1,i_2,\ldots ,i_N)$.
Given two such reduced expressions $\mathbf{i}$ and $\mathbf{i}'$,
we write $\mathbf{i}\sim
\mathbf{i}'$ if there is a sequence of commutations (of the form
$s_is_j=s_js_i$ with $|i-j|>1$) which, when applied to $\mathbf{i}$, give
$\mathbf{i}'$. This is an equivalence relation on the set of reduced
expressions for $w_0$, and the equivalence classes are called
commutation classes.

Let $Q=(Q_0,Q_1,s,t)$ be a quiver, i.e. a finite oriented graph with set of 
vertices $Q_0$, set of arrows $Q_1$, and maps $s,t:Q_1\rightarrow Q_0$
associating to each arrow its source and target, respectively.
If $Q$ is a quiver and $i$ is a sink in $Q$ (i.e. all the arrows incident with
$i$ have $i$ as target), we denote by $s_i(Q)$ the quiver with all of the
arrows incident with $i$ reversed. If $\mathbf{i}$ is a reduced expression
for $w_0$, we say (following Lusztig~\cite{lusztig2}) that $\mathbf{i}$ is
compatible with $Q$ if $i_1$ is a sink in $Q$, $i_2$ is a sink in
$s_{i_1}(Q)$, $i_3$ is a sink in $s_{i_2}s_{i_1}(Q)$, $\ldots$, $i_N$ is a
sink in $s_{i_{N-1}}s_{i_{N-2}}\cdots s_{i_1}(Q)$. It is known that if
$Q$ is any quiver of type $A_n$, then there is always at least one reduced
expression $\mathbf{i}$ compatible with $Q$, and that the set of reduced
expressions compatible with $Q$ is the commutation class of $\mathbf{i}$.

Berenstein, Fomin and Zelevinsky give a nice description of a reduced
expression compatible with a given quiver $Q$ in type $A_n$
in~\cite[\S4.4.3]{bfz1}.
Suppose that $Q$ is such a quiver. Number the edges of the quiver from $1$
to $n-1$, starting from the left hand end. Berenstein, Fomin and
Zelevinsky construct an arrangement as follows.
Consider a square in the plane, with horizontal
and vertical sides. We will draw $n+1$ `pseudo-lines' in this square.
Put $n+1$ points onto the left-hand edge of the square, equally spaced,
numbered $1$ to $n+1$ from top to bottom, so that $1$ and $n+1$ are at
the corners. Do the same for the right-hand edge, but number the points
from bottom to top. ${\rm Line}_h$ will join point $h$ on the left with
point $h$ on the right. For $h=1,n+1$, ${\rm Line}_h$ will be a diagonal
of the square. For $h\in [2,n]$, ${\rm Line}_h$ will be a union of two
line segments of slopes $\pi /4$ and $-\pi /4$. There are therefore precisely
two possibilities for $\mbox{Line}_h$. If edge $h-1$ in $Q$ is oriented to
the left, the left segment has positive slope, while the right one has
negative slope; if edge $h-1$ is oriented to the right, it goes the other
way round.

\begin{figure}[htbp]
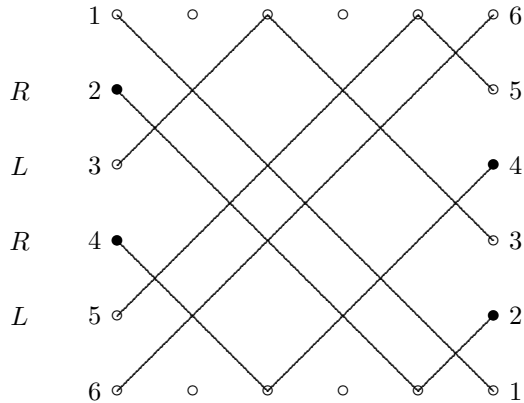

\beginpicture

\setcoordinatesystem units <1cm,1cm>             
\setplotarea x from -4 to 7, y from 0 to 7       
\linethickness=0.5pt           

\multiput{$\circ$} at 1 1 *5 0 1 /
\multiput{$\circ$} at 6 1 *5 0 1 /
\multiput{$\circ$} at 2 1 *3 1 0 /
\multiput{$\circ$} at 2 6 *3 1 0 /

\put{$\bullet$} at 1 3 %
\put{$\bullet$} at 1 5 %
\put{$\bullet$} at 6 2 %
\put{$\bullet$} at 6 4 %
\put{$1$}[c] at 0.7 6 %
\put{$2$}[c] at 0.7 5 %
\put{$3$}[c] at 0.7 4 %
\put{$4$}[c] at 0.7 3 %
\put{$5$}[c] at 0.7 2 %
\put{$6$}[c] at 0.7 1 %
\put{$6$}[c] at 6.3 6 %
\put{$5$}[c] at 6.3 5 %
\put{$4$}[c] at 6.3 4 %
\put{$3$}[c] at 6.3 3 %
\put{$2$}[c] at 6.3 2 %
\put{$1$}[c] at 6.3 1 %

\put{$R$}[c] at -0.3 5 %
\put{$L$}[c] at -0.3 4 %
\put{$R$}[c] at -0.3 3 %
\put{$L$}[c] at -0.3 2 %

\setlinear \plot 1 6    6 1 / %
\setlinear \plot 1 5    5 1 / %
\setlinear \plot 5 1    6 2 / %

\setlinear \plot 1 4    3 6 / %
\setlinear \plot 3 6    6 3 / %

\setlinear \plot 1 3    3 1 / %
\setlinear \plot 3 1    6 4 / %

\setlinear \plot 1 2    5 6 / %
\setlinear \plot 5 6    6 5 / %

\setlinear \plot 1 1    6 6 / %

\endpicture
\caption{The Berenstein, Fomin and Zelevinsky diagram for the quiver $RLRL$}
\label{bfzdiagram}
\end{figure}

{\bf Example:} Berenstein, Fomin and Zelevinsky give the following example.
Consider the case $A_5$, with $Q$ given by the quiver with
vertices $1,2,3,4,5$ and arrows $1\rightarrow 2$, $2\leftarrow 3$,
$3\rightarrow 4$, and $4\leftarrow 5$. (We shall denote such a quiver by
$RLRL$, where an $R$ (respectively, $L$) denotes an edge oriented to the
right (respectively, left).) The Berenstein-Fomin-Zelevinsky arrangement for
this quiver is shown in Figure~\ref{bfzdiagram}.

If $\mathbf{i}$ is a reduced expression for $w_0$ in type $A_n$, then the
{\em chamber diagram} for $\mathbf{i}$ is given by a set of pseudolines
numbered from $1$ to $n$. Two sets of points numbered $1$ to $n$ are
arranged on the vertical lines of a square as in the above picture.
Underneath the square, the simple reflections in $\mathbf{i}$ are written from
left to right. The $i$th pseudoline then links the point marked $i$ on the
left with the point marked $i$ on the right, in such a way that immediately
above the simple reflection $i_j$ from $\mathbf{i}$ pseudolines
$i_j$ and $i_j+1$ cross. See~\cite[1.4]{bfz1} for details. Berenstein,
Fomin and Zelevinsky prove the following result.

\begin{prop} \label{quiverarrangement} (Berenstein, Fomin and Zelevinsky)
A reduced expression $\bf i$ for $w_0$ is compatible with the quiver $Q$ if
and only if its chamber diagram is isotopic to the arrangement
defined above corresponding to $Q$.
\end{prop}

\noindent {\bf Proof:} see~\cite[\S4.4.3]{bfz1}.~$\square$

Since each simple reflection appearing in such a reduced expression
$\mathbf{i}$ corresponds to a crossing of two pseudo-lines in the diagram,
each pseudo-line in the diagram
(consisting of a part of positive slope and a part of negative slope) gives
rise to some of the simple reflections appearing in $\mathbf{i}$; each
simple reflection appears twice in this way as two lines crossing correspond
to a simple reflection. We can remove this duplication by counting, for each
line, only the simple reflections which arise during the part of the line
which is of positive slope. In this way, each edge of the quiver,
which corresponds to a line in the diagram, gives rise to some of the simple
reflections in the reduced expression $\mathbf{i}$; we include also the
line from the bottom left of the diagram to the top right.
Each simple reflection arises for a unique edge of the quiver (or comes from
the extra line).

Number the edges of $Q$ from $1$ to $n-1$, starting from the left. Suppose
that the edges $l_1,l_2,\ldots, l_a$ all point to the left, and that edges
$r_1,r_2,\ldots ,r_b$ all point to the right, and that every edge is one of
these, where $l_1<l_2<\cdots <l_a$ and $r_1<r_2<\cdots <r_b$.
For $m\in \mathbb{N}$, denote by $(m\searrow 1)$ the sequence
$m,m-1,\ldots ,2,1$. Then it is easy to see that the above construction
shows that the reduced expression
$$\mathbf{i}(Q)=(l_1\searrow 1)(l_2\searrow 1)\cdots (l_a\searrow 1)
(n\searrow 1)
(n\searrow n+1-r_b)(n\searrow n+1-r_{b-1})\cdots (n\searrow n+1-r_1)$$
is compatible with $Q$; it follows that the reduced expressions
compatible with $Q$ are precisely those commutation equivalent to
$\mathbf{i}(Q)$.

\section{Monomial Bases arising from Representations of Quivers}
\label{representations}

Let $Q=(Q_0,Q_1,s,t)$ be a quiver. We assume $Q$ to be of Dynkin type,
which means that the unoriented graph $\Delta$ underlying $Q$ is a disjoint 
union of Dynkin diagrams of type $A$, $D$, $E$.
Let $k$ be an arbitrary field. Then we can form the path algebra $kQ$ of $Q$ 
over $k$, which has the paths in $Q$ as a $k$-basis
(including an 'empty' path for each vertex $i\in Q_0$), and multiplication of 
paths given by concatenation if possible, and zero
otherwise. This is a finite dimensional $k$-algebra since $Q$, being of
Dynkin type, has no oriented cycles. We form the
category $\bmod kQ$ of finite-dimensional representations of $kQ$. The 
isoclasses of simple objects $S_i$ in $\bmod kQ$ correspond
bijectively to the vertices $i\in Q_0$ of $Q$. Let $\N Q$ be the free abelian 
semigroup spanned by elements $\alpha_i$ for $i\in Q_0$; it
can be identified with the positive root lattice of the root system $R$ of 
type $\Delta$. For a representation $M\in\bmod kQ$, we denote by
$d_i$ for $i\in Q_0$ the Jordan-H\"older multiplicity of the simple $S_i$ in 
$M$.  This allows us to define a map $\dimv$ from the set of
isoclasses in $\bmod kQ$ to $\N Q$ by $\dimv(M)=\sum_{i\in Q_0}d_i\alpha_i$. 
The fundamental result in the theory of Dynkin quivers is:

\begin{theorem} \label{gabriel} (Gabriel)\\
The map $\dimv$ induces a bijection between isoclasses of indecomposable 
objects $X_\alpha$ in $\bmod kQ$ and positive roots
$\alpha\in R^+\subset\N Q$ for type $\Delta$.
\end{theorem}

The Auslander-Reiten quiver of the algebra $kQ$ is defined as the oriented 
graph having the isoclasses of indecomposable representations
of $kQ$ as vertices, and arrows corresponding to irreducible maps (i.e. 
morphisms in $\bmod kQ$ between indecomposable objects which can not be
factored into a composition of non-split maps); see~\cite{ars} for details.
An explicit construction will be given in section \ref{lusztigcones}.

The following definition was first introduced in~\cite{reineke3}:

\begin{definition} \label{directedpartition}
A partition $R^+=I_1\cup\ldots\cup I_s$ of $R^+$ into disjoint subsets $I_k$ 
is called directed if\\
(i) $\Ext(X_\alpha,X_\beta)=0$ for all $\alpha,\beta$ in the same part $I_k$,\\
(ii) $\Ext(X_\alpha,X_\beta)=0$ and $\Hom(X_\beta,X_\alpha)=0$ if
$\alpha\in I_k$, $\beta\in I_l$, where $1\leq k<l\leq s$.
\end{definition}

The existence of (several!) directed partitions of a given quiver $Q$ can be 
seen using Auslander-Reiten theory (see~\cite{ars}):
if $\Hom(U,V)\not=0$ (resp. $\Ext(V,U)\not=0$) for indecomposables 
$U,V\in\bmod kQ$, then there exists a path (resp. a proper path)
from $[U]$ to $[V]$ in
the Auslander-Reiten quiver of $kQ$. But since this graph is directed, we can 
enumerate the isoclasses of indecomposables in
$\bmod kQ$ as $[U_1],\ldots,[U_\nu]$, such that $\Hom(U_p,U_q)=0$ for $p>q$, 
and $\Ext(U_p,U_q)=0$ for $p\leq q$. Define roots $\alpha^p$ by 
$[U_p]=[X_{\alpha^p}]$. By definition, the partition
$R^+=\{\alpha^1\}\cup\ldots\cup\{\alpha^\nu\}$ is directed. All other
directed  partitions can be constructed by coarsening such a partition.

Fix a directed partition $R^+=I_1\cup\ldots\cup I_s$ from now on. We will 
associate to it a sequence $\mathbf{i}=i_1\ldots i_t$, as well as
a function $D$ from the set of isoclasses in $\bmod kQ$ to $\N^t$. 
Enumerate the vertices $Q_0$ of $Q$ as $Q_0=\{1\ldots n\}$ such that
the existence of an arrow $i\rightarrow j$ in $Q$ implies $i<j$. For each 
$p=1\ldots s$, write the 
subset of $Q_0$ consisting of all vertices $i$ such that $\alpha_i$
appears with non-zero coefficient in some root $\alpha\in I_p$ as
$\{i^p_1,\ldots,i^p_{t_p}\}$, increasing with respect to the above
defined ordering on $Q_0$. Then the sequence $\mathbf{i}$ is defined as
$$\mathbf{i}=i^1_1\ldots i^1_{t_1}\, i^2_1\ldots i^2_{t_2}\,\ldots\, 
i^s_1\ldots i^s_{t_s}.$$
Given an isoclass $[M]$ in $\bmod kQ$, we can write $M$ as 
$M=\oplus_{\alpha\in R^+}X_\alpha^{c_\alpha}$, using Gabriel's Theorem and
Krull-Schmidt. Write $s\in \alpha$ to indicate that the simple root $\alpha_s$
appears in the expression for $\alpha$ as a sum of simple roots.
Let $D(M)$ be the tuple $\mathbf{a}=(a_j)_{j=1,2,\ldots ,k}$
defined as follows. If $\alpha^j$ lies in $I_p$, then
\begin{equation}\label{Ddef}
a_j =\sum_{\alpha\in I_p,\ i_j\in \alpha}c_{\alpha}.
\end{equation}
This defines a function $D$ from the set of isoclasses in $\bmod kQ$ 
to $\N^N$, which is obviously additive, i.e.
$D(M\oplus N)=D(M)+D(N)$.
Identifying the set of isoclasses in $\bmod kQ$ with $\N R^+$ via 
$[\oplus_\alpha X_\alpha^{c_\alpha}]\mapsto \sum_{\alpha}c_\alpha\alpha$,
we thus get a linear function $D:\Z R^+\rightarrow \Z^t$. The main 
result of this paper is that, for suitable choices of $Q$ and the
directed partition, the function $D$ coincides with one of the 
reparametrization functions
$(S_\mathbf{i}^\mathbf{j})^{-1}$ on an explicitly described region.

The original use of the function $D$ lies in the following theorem 
(see~\cite{reineke3}):

\begin{theorem}\label{monomialbases} (Reineke)
Writing $\mathbf{i}=i_1\ldots i_t$ and $D=(D_1,\ldots,D_t)$, the set
$$\{F_{i_1}^{D_1(\mathbf{c})}\cdot\ldots F_{i_t}^{D_t(\mathbf{c})}\in U^-
\, :\, \mathbf{c}\in \N R^+\}$$
is a basis for $U^-$.
\end{theorem}

Moreover, these monomial bases for $U^-$ have good properties with respect to 
base change: the base change coefficients to a PBW basis (resp. to the 
canonical basis) form upper unitriangular matrices with respect to a certain 
ordering (namely, related to the degeneration ordering on quiver representations), and these coefficients have representation theoretic (resp. 
geometric) interpretations (see~\cite{reineke3},~\cite{reineke4}).

The function $D$ is not invertible in general, but it is for special directed 
partitions, called regular in~\cite{reineke3}. We now
introduce a certain class of directed partitions for quivers of type $A$,
and compute the inverse of $D$ in these special cases.

\begin{figure}[htbp]
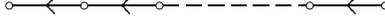


\beginpicture

\setcoordinatesystem units <1cm,1cm>             
\setplotarea x from -2.5 to 12, y from 1.7 to 2.3       

\scriptsize{

\multiput {$\circ$} at 3   2 *2 1 0 /      %
\multiput {$\circ$} at 7   2 *1 1 0 /      %

\putrule from 3.05 2 to 3.95 2  %
\putrule from 4.05 2 to 4.95 2  %
\putrule from 7.05 2 to 7.95 2  %

\setlinear \plot  3.6 2.1  3.5 2  / %
\setlinear \plot  3.6 1.9  3.5 2  / %
\setlinear \plot  4.6 2.1  4.5 2  / %
\setlinear \plot  4.6 1.9  4.5 2  / %
\setlinear \plot  7.6 2.1  7.5 2  / %
\setlinear \plot  7.6 1.9  7.5 2  / %

\setdashes <2mm,1mm>          %

\putrule from 5.05 2 to 6.95 2 %
}
\endpicture
\caption{The quiver compatible with $\mathbf{k}$}
\label{linearquiver}
\end{figure}

We thus assume that $Q$ is of type $A_n$, and we
fix once and for all the reduced expression for $w_0$ given by
$\mathbf{k}=(1,2,1,3,2,1,\ldots ,n,n-1,\ldots ,1)$. This is one of the most
regular reduced expressions, and is compatible with the quiver
$Q_{\mathbf{k}}$ given in Figure~\ref{linearquiver}. We conjecture that our
results hold, with similar proofs, if this quiver is replaced by an arbitrary
Dynkin quiver of type $A_n$, but the Kashiwara operators are in this general
case more difficult to describe combinatorially. In the special case we are
considering, it is possible to handle this combinatorics.
See~\cite{reineke1} for a description
of the combinatorics of Kashiwara operators on canonical basis elements using
parametrizations of the canonical basis arising from PBW-bases of
$U^-$ corresponding to reduced expressions compatible with arbitrary quivers.

We can construct the Auslander-Reiten quiver of our quiver $Q$ in the
following way (see~\cite{ars}). Let
${\mathbb{N}}Q$ be the quiver with
vertices $\mathbb{N}\times \{1,2,\ldots ,n\}$. Whenever there is an arrow
$i\rightarrow j$ in $Q$, we draw one arrow $(z,i)\rightarrow (z,j)$ and
one arrow $(z,j)\rightarrow (z+1,i)$ for each $z\in \mathbb{N}$. Define
$A(Q)$ to be the full subquiver of $\mathbb{N}Q$ consisting of all vertices
$(z,i)$ such that $1\leq z\leq (h+a_i-b_i)/2$ where, for each $i\in
\{1,2,\ldots ,n\}$, $a_i$ (respectively $b_i$), is the number of arrows in the
unoriented path from $i$ to $\sigma(i)$ that are directed towards $i$
(respectively $\sigma(i)$). Here, $\sigma$ is the unique permutation of
the vertices of $Q$ such that $w_0(\alpha_i)=-\alpha(\sigma(i))$, and $h$
is the Coxeter number. Then $A(Q)$ is the Auslander-Reiten quiver of $Q$.

A reduced expression $\mathbf{i}$ for $w_0$ defines an ordering on the set
$\Phi^+$ of positive roots of the root system associated to $W$.
We write $\alpha^j=s_{i_1}s_{i_2}\cdots s_{i_{j-1}}(\alpha_{i_j})$ for
$j=1,2,\ldots ,N$. Then $\Phi^+=\{\alpha^1,\alpha^2,\ldots \alpha^N\}$.
For ${{\mathbf{c}}}=(c_1,c_2,\ldots,c_N)\in \mathbb{Z}^N$, write
$c_{\alpha^j}=c_j$.
If $\alpha=\alpha_{ij}:=\alpha_i+\alpha_{i+1}+\cdots +\alpha_{j}$ with $i<j$,
we also write $c_{ij}$ for $c_{\alpha_{ij}}$.

\begin{figure}[htbp]
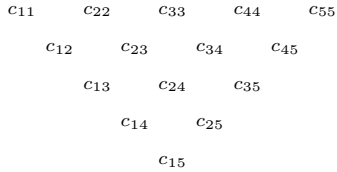


\beginpicture

\setcoordinatesystem units <0.5cm,0.5cm>             
\setplotarea x from -15 to 5, y from -0.3 to 4.3       

\scriptsize{

\put {$c_{11}$} at -4 4       %
\put {$c_{22}$} at -2 4      %
\put {$c_{33}$} at 0 4      %
\put {$c_{44}$} at 2 4      %
\put {$c_{55}$} at 4 4      %
\put {$c_{12}$} at -3 3      %
\put {$c_{23}$} at -1 3      %
\put {$c_{34}$} at 1 3      %
\put {$c_{45}$} at 3 3       %
\put {$c_{13}$} at -2 2      %
\put {$c_{24}$} at 0 2      %
\put {$c_{35}$} at 2 2      %
\put {$c_{14}$} at -1 1      %
\put {$c_{25}$} at 1 1      %
\put {$c_{15}$} at 0 0      %
}
\endpicture
\caption{Array of elements of $\mathbb{N}^{15}$}
\label{triangle}
\end{figure}

Let $\alpha^1,\alpha^2,\ldots ,\alpha^N$ be the ordering induced on
$\Phi^+$ by $\mathbf{k}$. We can write an element of $\mathbb{N}^N$ as
$c=(c_{ij})$, where $1\leq i\leq j\leq n$, using the above.
The canonical basis is parametrized, via the Lusztig parametrization
$\phi_{\mathbf{k}}:\mathbf{B} \rightarrow \mathbb{N}^N$. We can write an
element of $\N^N$ as an array based on the Auslander-Reiten quiver for
$Q_{\mathbf{k}}$. We write $c_{ij}$ in place of the module corresponding to
the positive root $\alpha_i+\alpha_{i+1}+\cdots \alpha_{j}$,
with $c_{11},c_{22},\ldots ,c_{nn}$ on the first row,
$c_{12},c_{23},\ldots ,c_{n-1,n}$ along the second, interspersing the first
row elements, and so on, until $c_{1,n}$ on the last row. For example,
when $n=5$, see Figure~\ref{triangle}.
Suppose that $Q'$ is another quiver of type $A_n$. We can define a directed
partition of $A(Q)$ in the following way. Fix $z\in \mathbb{N}$, and set
$z_1=z$. Let $v_1:=(z_1,1)$ be a vertex of
$\mathbb{N}Q$. For $i=2\ldots n$, define a vertex $v_i=(z_i,i)$ of
$\mathbb{N}Q$ inductively, as follows. If $i-1\rightarrow i$ is an arrow in
$Q'$, then let $v_i$ be the head of the unique arrow with source $v_{i-1}$.
If $i-1\leftarrow i$ is an arrow in $Q'$, then let $v_i$ be the source of the
unique arrow with head $v_{i-1}$.
Let $S_z=(v_1,v_2,\ldots ,v_n)\subseteq \mathbb{N}Q$. Then it is clear that
$\mathbb{N}Q=\cup_{z\in \mathbb{N}}S_z$. It follows that
$A(Q)=\cup_{z\in \mathbb{N}}(S_z\cap A(Q))$; note that this decomposition must
be finite. It is clear from the construction that this is a directed
partition of $A(Q)$. Let $T_z=S_z\cap A(Q)$. We call each $S_z$ a {\em slice}
of $\mathbb{N}Q$. If $T_z$ is non-empty, we call it a slice of $A(Q)$.

\begin{figure}[htbp]
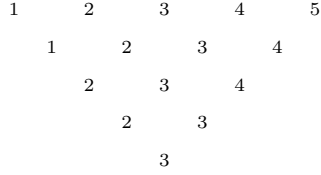


\beginpicture

\setcoordinatesystem units <0.5cm,0.5cm>             
\setplotarea x from -15 to 5, y from -0.3 to 4.3       

\scriptsize{

\put {$1$} at -4 4       %
\put {$2$} at -2 4      %
\put {$3$} at 0 4      %
\put {$4$} at 2 4      %
\put {$5$} at 4 4      %
\put {$1$} at -3 3      %
\put {$2$} at -1 3      %
\put {$3$} at 1 3      %
\put {$4$} at 3 3       %
\put {$2$} at -2 2      %
\put {$3$} at 0 2      %
\put {$4$} at 2 2      %
\put {$2$} at -1 1      %
\put {$3$} at 1 1      %
\put {$3$} at 0 0      %
}
\endpicture
\caption{Slice structure of $A(Q_{\mathbf{k}})$ corresponding to the quiver
$RLRL$}
\label{sliceRLRL}
\end{figure}

{\bf Example:} Consider the case $A_5$, with $Q'=RLRL$.
The corresponding decomposition of $A(Q)$ into slices is given in
Figure~\ref{sliceRLRL}. Each vertex of $A(Q)$ is denoted by a number,
indicating the number of the slice it lies in.

Let us now specialise to the case where $Q=Q_{\mathbf{k}}$ is the quiver
with arrows $i \leftarrow i+1$ for $i=1,2,\ldots ,n-1$, compatible with
$\mathbf{k}$. Let $Q'$ be an arbitrary quiver.
Note that we can identify $(i,j)$, for $1\leq i\leq j\leq n$, with the
indecomposable module with dimension vector $\alpha_i+\cdots +\alpha_j$,
so each such pair lies in a corresponding slice $T_z$.
If $v=(z,a)\in A(Q)$, let $i(v),j(v)$ be
such that the dimension vector of the module at vertex $v$ is
$\alpha_{i(v)}+\cdots +\alpha_{j(v)}$.

Recall the function $D:\Z R^+\rightarrow \Z^k$ associated in section 
\ref{representations} to a quiver and a directed partition. Consider the
special case where $Q=Q_{\mathbf{k}}$, and where the directed partition
is associated  to another quiver $Q'$ of type $A$ as above.
In this case it is easy to see that the sequence $\mathbf{i}$ is given by
$$\mathbf{i}(Q')=(l_1\searrow 1)(l_2\searrow 1)\cdots (l_a\searrow 1)
(n\searrow 1)
(n\searrow n+1-r_b)(n\searrow n+1-r_{b-1})\cdots (n\searrow n+1-r_1),$$
a reduced expression for the longest word in the Weyl group compatible with
$Q'$.
Here, each of the bracketed parts of $\mathbf{i}(Q')$ arises from a
part of the directed partition associated to $Q'$.

In this special case, we can easily prove the invertibility of $D$.

\begin{lemma} \label{Dinvertibility}
Suppose that $Q=Q_{\mathbf{k}}$, and that $D$ is the function associated to
the directed partition of $A(Q)$ corresponding to another quiver $Q'$ of
type $A$. Then the function $D$ is invertible with inverse function
$E=D^{-1}$. Moreover, the components $c_\alpha$ of
$\mathbf{c}=(c_\alpha)_\alpha=E(\mathbf{a})$ for some 
$\mathbf{a}\in\N^N$ and $\alpha\in R^+$ are of the form
$c_\alpha=a_k-a_l$ or $c_\alpha=a_k$.
\end{lemma}

\noindent {\bf Proof:} Recall from Section~\ref{representations} that we can 
write
$$\mathbf{i}=(i^1_1\ldots i^1_{t_1}\,\ldots\, i^s_1\ldots i^s_{t_s})$$
and that $D(\mathbf{c})$ is the tuple $\mathbf{a}=(a_j)_{j=1,2,\ldots ,k}$
such that if the root $\alpha^j$ lies in $I_p$, then
$$a_j =\sum_{\alpha\in I_p,\ i_j\in \alpha}c_{\alpha}.$$
In our case, it is easily seen from the definition of the directed partition 
that we can write
$$I_p=\{\alpha^p_1,\ldots,\alpha^p_{t_p'}\},$$
such that the root $\alpha^p_u$ is of length $u$, and 
$\alpha^p_u=\alpha^p_{u-1}+\alpha_{j_u}$ for some simple root $\alpha_{j_u}$ 
(we formally set $\alpha^p_0=0$). It follows that 
$\alpha^p_u=\alpha_{j_1}+\ldots+\alpha_{j_u}$, which in particular implies 
$t_p'=t_p$. We also see that
$$\{i^p_1,\ldots,i^p_{t_p}\}=\{j_1,\ldots,j_{t_p}\};$$
denote by $\sigma$ the permutation defined by $j_{\sigma(u)}=i^p_u$. Then we 
can compute $a^p_u$ as:
$$a^p_u=\sum_{\alpha\in I_p,\ i_{k_u}\in \alpha}c_{\alpha}=
c_{\alpha^p_{\sigma u}}+\ldots+c_{\alpha^p_{t_p}}.$$
It follows that
$$c_{\alpha^p_u}=a^p_{\sigma^{-1}u}-a^p_{\sigma^{-1}(u+1)}\mbox{ if 
$u\not=t_p$, and }c_{\alpha^p_u}=a^p_{\sigma^{-1}u}\mbox{ otherwise}.$$
This proves the claimed properties of $D$.

{\bf Remark on Notation}

In the sequel, we shall consider directed partitions in $A(Q_{\mathbf{k}})$
(so that, in the above, $Q=Q_{\mathbf{k}}$), arising from an arbitrary
quiver $Q$ of type $A_n$ (which was denoted $Q'$ above).

\section{The Lusztig cones} \label{lusztigcones}

Lusztig~\cite{lusztig7}
introduced certain regions which, in low rank, give rise to
canonical basis elements of a particularly simple form.
The {\em Lusztig cone} corresponding
to a reduced expression $\mathbf{i}$ for $w_0$ is defined to be the set of
points ${{\mathbf{a}}}\in {\N}^N$
satisfying the following inequalities: \\
(*) For every pair $s,s'\in [1,k]$ with $s<s'$, $i_s=i_{s'}=i$ and
$i_p\not=i$ whenever $s<p<s'$, we have
$$(\sum_p a_p)-a_s - a_{s'}\geq 0,$$
where the sum is over all $p$ with $s<p<s'$ such that $i_p$ is joined to $i$
in the Dynkin diagram. We shall denote this cone by $L_{st}(\mathbf{i})$
(as we shall regard it as a set of strings of ${\B}$ in direction
$\mathbf{i}$).

It was shown by Lusztig~\cite{lusztig7} that, in type $A_n$, if $\mathbf{a}\in
L_{st}(\mathbf{i})$ then the monomial $F_{i_1}^{(a_1)}F_{i_2}^{(a_2)}\cdots
F_{i_N}^{(a_N)}$ lies in the canonical basis ${\B}$, provided $n=1,2,3$.
The first author~\cite{marsh7} showed that this remains true if $n=4$, but it
is false for $n\geq 5$ by~\cite{marsh7},~\cite{reineke2}.
The Lusztig cones have been studied in the
papers~\cite{cartermarsh},~\cite{marsh8} and~\cite{marsh9} in type $A$ for
every reduced expression
$\bf i$ for the longest word, and have also been studied by Bedard
in~\cite{bedard2} for arbitrary finite
(simply-laced) type for reduced expressions compatible with a quiver whose
underlying graph is the Dynkin diagram. Bedard describes these vectors
using the Auslander-Reiten quiver of the quiver and homological algebra,
showing they are closely connected to the representation theory of the quiver.

\section{The Degeneration Cones} \label{mrcones}

We define the cone $C_{PBW}(Q)\subseteq \mathbb{N}^N$
corresponding to a quiver $Q$ to be the set of points
$c=(c_{ij})$ satisfying the inequalities (C1) and (C2) below.
We define the {\em degeneration cone} corresponding to $Q$ to be
the cone $C_{st}(Q)=D(C_{PBW}(Q))$.

Define a {\em component} of $Q$ to be a maximal full subgraph $X$ of
$Q$ subject to the condition that all of the arrows of $X$ point in the
same direction. Call $X$ a {\em left} (respectively, {\em right}) component
of $Q$ if its arrows all point to the left (respectively, right).
If $X$ is a component of $Q$ (left or right), let $S_z(X)$ (respectively,
$T_z(X)$) denote the part of the slice $S_z$ (respectively, $T_z$)
corresponding to $X$.

\begin{figure}[htbp]
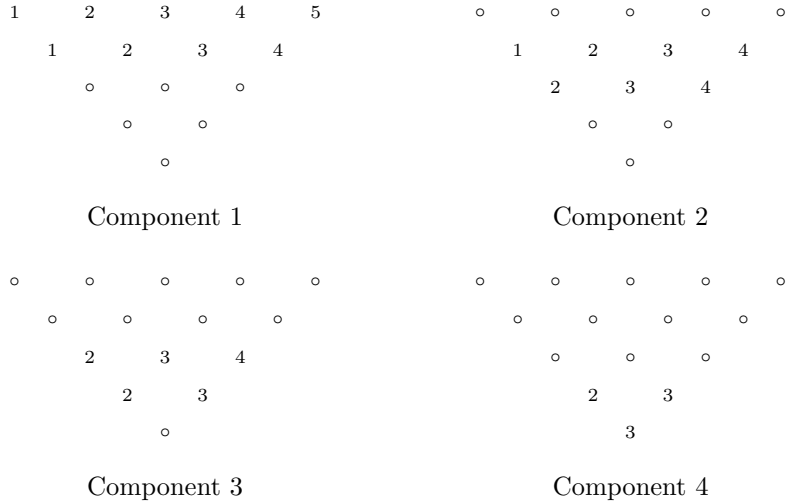

\begin{minipage}{2in}

\beginpicture

\setcoordinatesystem units <0.5cm,0.5cm>             
\setplotarea x from -8 to 4, y from -0.3 to 4.3       

\scriptsize

\put {$1$} at -4 4       %
\put {$2$} at -2 4      %
\put {$3$} at 0 4      %
\put {$4$} at 2 4      %
\put {$5$} at 4 4      %
\put {$1$} at -3 3      %
\put {$2$} at -1 3      %
\put {$3$} at 1 3      %
\put {$4$} at 3 3       %
\put {$\circ$} at -2 2      %
\put {$\circ$} at 0 2      %
\put {$\circ$} at 2 2      %
\put {$\circ$} at -1 1      %
\put {$\circ$} at 1 1      %
\put {$\circ$} at 0 0      %
\normalsize
\put {Component $1$} at 0 -1.5 %
\endpicture
\end{minipage}
\begin{minipage}{2in}

\beginpicture

\setcoordinatesystem units <0.5cm,0.5cm>             
\setplotarea x from -8 to 4, y from -0.3 to 4.3       

\scriptsize
\put {$\circ$} at -4 4       %
\put {$\circ$} at -2 4      %
\put {$\circ$} at 0 4      %
\put {$\circ$} at 2 4      %
\put {$\circ$} at 4 4      %
\put {$1$} at -3 3      %
\put {$2$} at -1 3      %
\put {$3$} at 1 3      %
\put {$4$} at 3 3       %
\put {$2$} at -2 2      %
\put {$3$} at 0 2      %
\put {$4$} at 2 2      %
\put {$\circ$} at -1 1      %
\put {$\circ$} at 1 1      %
\put {$\circ$} at 0 0      %
\normalsize
\put {Component $2$} at 0 -1.5 %
\endpicture
\end{minipage}

\vspace*{0.5cm}
\begin{minipage}{2in}

\beginpicture

\setcoordinatesystem units <0.5cm,0.5cm>             
\setplotarea x from -8 to 4, y from -0.3 to 4.3       

\scriptsize

\put {$\circ$} at -4 4       %
\put {$\circ$} at -2 4      %
\put {$\circ$} at 0 4      %
\put {$\circ$} at 2 4      %
\put {$\circ$} at 4 4      %
\put {$\circ$} at -3 3      %
\put {$\circ$} at -1 3      %
\put {$\circ$} at 1 3      %
\put {$\circ$} at 3 3       %
\put {$2$} at -2 2      %
\put {$3$} at 0 2      %
\put {$4$} at 2 2      %
\put {$2$} at -1 1      %
\put {$3$} at 1 1      %
\put {$\circ$} at 0 0      %
\normalsize
\put {Component $3$} at 0 -1.5 %
\endpicture
\end{minipage}
\begin{minipage}{2in}

\beginpicture

\setcoordinatesystem units <0.5cm,0.5cm>             
\setplotarea x from -8 to 4, y from -0.3 to 4.3       

\scriptsize
\put {$\circ$} at -4 4       %
\put {$\circ$} at -2 4      %
\put {$\circ$} at 0 4      %
\put {$\circ$} at 2 4      %
\put {$\circ$} at 4 4      %
\put {$\circ$} at -3 3      %
\put {$\circ$} at -1 3      %
\put {$\circ$} at 1 3      %
\put {$\circ$} at 3 3       %
\put {$\circ$} at -2 2      %
\put {$\circ$} at 0 2      %
\put {$\circ$} at 2 2      %
\put {$2$} at -1 1      %
\put {$3$} at 1 1      %
\put {$3$} at 0 0      %
\normalsize
\put {Component $4$} at 0 -1.5 %
\endpicture
\end{minipage}
\label{subsliceRLRL}
\caption{The subsets $T_z(X)$ of $A(Q_{\mathbf{k}})$ corresponding to
components of the quiver $RLRL$}
\end{figure}

{\bf Example: } Consider the example with $Q=RLRL$ given above. Then $Q$ has
$4$ components, each containing one edge. Two are right components, and
two are left components. For each component $X$, we indicate the subsets of
slices, $T_z(X)$, in $A(Q)$; see Figure~\ref{subsliceRLRL}. In each case,
the numbers $z$ denote elements of the sets $T_z(X)$, and the empty
circles denote elements not in any subset $T_z(X)$.

Firstly, let $X$ be a left component of $Q$, and suppose $T_z,T_{z+1}$
are consecutive slices in $A(Q)$ such that $S_z(X)=T_z(X)$ and $S_{z+1}(X)=
T_{z+1}(X)$ (i.e. the sets $S_z(X)$ and $S_{z+1}(X)$ of $\mathbb{N}Q$ are
contained entirely inside $A(Q)$).
We can order $T_z(X)$ linearly by the second entry
of the vertices appearing in it (recall from Section~\ref{representations}
that each vertex is regarded
as a pair $(r,i)$, where $r\in \mathbb{N}$ and $i$ is a vertex of $Q$, i.e.
an element of $\{1,2,\ldots ,n\}$).
Write $T_z(X)=\{x_1,x_2,\ldots x_k\}$.
Similarly, we can write $T_{z+1}(X)=\{y_1,y_2,\ldots ,y_k\}$.
For $t=1\ldots k$, let $i_t=i(x_t)$, and let $j_t=j(x_t)$. Note that
$i(y_t)=i_t+1$ and $j(y_t)=j_t+1$.
Then inequalities (C1) are:

(C1) $\sum_{r=a}^k c_{i_r,j_r}\geq \sum_{r=a}^k c_{i_r+1,j_r+1}$ for
$a=1\ldots k$.

Secondly, let $X$ be a right component of $Q$, and suppose $T_z,T_{z+1}$
are consecutive slices in $A(Q)$. We can order $T_z(X)$ linearly by the
second component of the vertices appearing in it; write
$T_z(X)=\{x_1,x_2,\ldots x_k\}$. Similarly, we can write
$T_{z+1}(X)=\{y_1,y_2,\ldots ,y_l\}$, where $k\geq l$, since $X$ is a right
component.
For $t=1\ldots k$, let $i_t=i(x_t)$, and let $j_t=j(x_t)$. Note that
$i(y_t)=i_t+1$ and $j(y_t)=j_t+1$ for $t=1\ldots l$.
Then inequalities (C2) are:

(C2) $c_{i_r,j_r}\geq c_{i_r+1,j_r+1}$ for $r=1\ldots l-1$.

The cone $C_{PBW}(Q)$ is defined to be the set of $(c_{ij})\in \mathbb{N}^N$
satisfying all of the inequalities (C1) and (C2).

{\bf Example: } In our running example (see the start of this section),
the inequalities (C1) and (C2) are given as follows.

Component $1$: $c_{11}\geq c_{22}\geq c_{33}\geq c_{44}\geq c_{55}$. \\
Component $2$: $c_{13}\geq c_{24}\geq c_{35}$ and $c_{13}+c_{23}\geq
c_{24}+c_{34}\geq c_{35}+c_{45}$. \\
Components $3$ and $4$ give rise to no inequalities in this example.

\section{The PBW-version of the Lusztig cones}

We now define a second cone $L_{PBW}(Q)$. It will be seen later that this cone
is the image, under the reparametrization map
$S_{\mathbf{i}(Q)}^{\mathbf{j}}$, of the Lusztig cone $L_{st}(Q)$. We'll
then see that $C_{PBW}(Q)$ (which will be seen to be the PBW-version of the
degeneration cone) contains the PBW-version of the Lusztig cone $L_{PBW}(Q)$.

Firstly, let $X$ be a left component of $Q$, and suppose $T_z,T_{z+1}$
are consecutive slices in $A(Q)$.
We can order $T_z(X)$ linearly by the second entry
of the vertices appearing in it; write $T_z(X)=\{x_1,x_2,\ldots x_k\}$.
Similarly, we can write $T_{z+1}(X)=\{y_1,y_2,\ldots ,y_l\}$, where $k\leq l$
(since $X$ is a left component).
For $t=1\ldots k$, let $i_t=i(x_t)$, and let $j_t=j(x_t)$. Note that
$i(y_t)=i_t+1$ and $j(y_t)=j_t+1$ for $1\leq t\leq k$.
Then the defining inequalities of
$L_{PBW}(Q)$ are

(L1) If $S_z(X)=T_z(X)$ and $S_{z+1}(X)=T_{z+1}(X)$ then
$\sum_{r=1}^k c_{i_r,j_r}\geq \sum_{r=1}^k c_{i_r+1,j_r+1}$, and \\
(L2) $c_{i_r,j_r}\leq c_{i_r+1,j_r+1}$ for $r=2\ldots k-1$.

Secondly, let $X$ be a right component of $Q$, and suppose $T_z,T_{z+1}$
are consecutive slices in $A(Q)$.
We can order $T_z(X)$ linearly by the second entry
of the vertices appearing in it; write $T_z(X)=(x_1,x_2,\ldots x_k)$.
Similarly, we can write $T_{z+1}(X)=(y_1,y_2,\ldots ,y_l)$, where $k\geq l$
(since $X$ is a right component).
For $t=1\ldots l$, let $i_t=i(x_t)$, and let $j_t=j(x_t)$. Note that
$i(y_t)=i_t+1$ and $j(y_t)=j_t+1$ for $1\leq t\leq l$.
Then the defining inequalities of
$L_{PBW}(Q)$ are

(L3) If $S_z(X)=T_z(X)$ and $S_{z+1}(X)=T_{z+1}(X)$ then
$\sum_{r=1}^k c_{i_r,j_r}\leq \sum_{r=1}^k c_{i_r+1,j_r+1}$, and \\
(L4) $c_{i_r,j_r}\geq c_{i_r+1,j_r+1}$ for $r=2\ldots l-1$.

Finally, if $X$ is a left component and is the leftmost component of $Q$,
then we have:

(L5) $c_{i_1,j_1}\leq c_{i_1+1,j_1+1}$,

and if $X$ is a right component and is the leftmost component of $Q$,
then we have:

(L6) $c_{i_1,j_1}\geq c_{i_1+1,j_1+1}$.

\section{The relationship between the Lusztig cones and the degeneration
cones in the PBW-parametrization}

We show in this section, that, for an arbitrary quiver $Q$ of type $A_n$,
we have $L_{PBW}(Q)\subseteq C_{PBW}(Q)$. Firstly, we show that certain
inequalities hold in $L_{PBW}$. Note that we will show in the next sections
that $L_{PBW}(Q)$ is the set of PBW-parameters of the Lusztig cone $L_{st}(Q)$

\begin{lemma} \label{inequalities}
Suppose $\mathbf{c}\in L_{PBW}(Q)$. Then the following inequalities hold:

If $X$ is a left component of $Q$, then (with the notation above),

(I) $c_{i_1,j_1}\leq c_{i_1+1,j_1+1}$, and \\
(II) $c_{i_k,j_k}\geq c_{i_k+1,j_k+1}$.

If $X$ is a right component, then

(III) $c_{i_1,j_1}\geq c_{i_1+1,j_1+1}$, and \\
(IV) $c_{i_k,j_k}\leq c_{i_k+1,j_k+1}$.

The inequalities (III) and (IV) only hold if they make sense --- i.e. if
$k=l$ (which is the case when $S_z(X)=T_z(X)$ and $S_{z+1}(X)=T_{z+1}(X)$).
\end{lemma}

{\bf Proof:}
Let $\mathbf{c}\in L_{PBW}(Q)$. Suppose first that the leftmost component
$X_1$ of $Q$ is a left component.
The inequality (I) for $X_1$ is the defining inequality (L5) of $L_{PBW}(Q)$.
Suppose that $S_z(X)=T_z(X)$ and $S_{z+1}(X)=T_{z+1}(X)$.
Then we have, from (L1) and (L2), that:

(a) $\sum_{r=1}^k c_{i_r,j_r}\geq \sum_{r=1}^k c_{i_r+1,j_r+1}$, and \\
(b) $c_{i_r,j_r}\leq c_{i_r+1,j_r+1}$ for $r=2\ldots k-1$.

Since (I) and (b) hold, if (II) were false, (a) would be false --- a
contradiction. Therefore (II) holds for $X_1$.

A similar argument shows that the inequalities (III) and (IV) hold for $X_1$
in the case when $X_1$ is a right component.

Suppose we are still in the case where $X_1$ is a left component.
Let $X_2$ be the component immediately to the right of $X_1$ (if it exists);
this must be a right component.
Note that the inequality (II) for $X_1$ is inequality (III) for $X_2$. We can
then argue as above for $X_2$ to deduce (IV) for $X_2$.

We can argue similarly in the case where $X_1$ is a right
component. In this way we can deduce the relevant inequalities for all
components of $Q$, by induction on the number of the component, starting from
the left, and we are done.~$\Box$

\begin{prop}
Let $Q$ be an arbitrary quiver of type $A_n$, and let $L_{PBW}(Q)$ and
$C_{PBW}(Q)$ be the cones as above. Then $L_{PBW}(Q)\subseteq C_{PBW}(Q)$.
\end{prop}

{\bf Proof:}
Let $\mathbf{c}\in L_{PBW}(Q)$.
Suppose that $X$ is any left component of $Q$, and that
$S_z(X)=T_z(X)$ and $S_{z+1}(X)=T_{z+1}(X)$. The following inequalities
hold:

(c) $\sum_{r=1}^k c_{i_r,j_r}\geq \sum_{r=1}^k c_{i_r+1,j_r+1}$, and \\
(d) $c_{i_r,j_r}\leq c_{i_r+1,j_r+1}$ for $r=1\ldots k-1$ and
$c_{i_k,j_k}\geq c_{i_k+1,j_k+1}$.

Inequality (c) comes from (L1), and inequalities (d) come from (L2) and
Lemma~\ref{inequalities}(II). We show that the inequalities (C1) all hold.
We already have the inequality $c_{i_k,j_k}\geq c_{i_k+1,j_k+1}$, from (d).
If we also had $c_{i_{k-1},j_{k-1}}+c_{i_k,j_k}\leq c_{i_{k-1}+1,j_{k-1}+1}+
c_{i_k+1,j_k+1}$, then this, together with the inequalities
$c_{i_r,j_r}\leq c_{i_r+1,j_r+1}$ for $r=1,\ldots ,k-2$ would give
$\sum_{r=1}^k c_{i_r,j_r}\leq \sum_{r=1}^k c_{i_r+1,j_r+1}$, contradicting
(c). Similarly, if we had
$\sum_{r=a}^k c_{i_r,j_r}\leq \sum_{r=a}^k c_{i_r+1,j_r+1}$, for some $1\leq
a\leq k-2$, the inequalities in (d) would give us
$\sum_{r=1}^k c_{i_r,j_r}\leq \sum_{r=a}^k c_{i_r+1,j_r+1}$, contradicting
(c). Hence
$\sum_{r=a}^k c_{i_r,j_r}\geq \sum_{r=a}^k c_{i_r+1,j_r+1}$, for any
$1\leq a\leq k-1$, and we see that the defining inequalities (C1) of
$C_{PBW}(Q)$ are satisfied.

Suppose next that $X$ is any right component of $Q$.
Then the following inequalities hold, using (L4) and
Lemma~\ref{inequalities}(III).

(e) $c_{i_r,j_r}\geq c_{i_r+1,j_r+1}$ for $r=1\ldots l-1$.

In this case, (e) contains all of the defining inequalities (C2)
of $C_{PBW}(Q)$. We thus see that $L_{PBW}(Q)\subseteq C_{PBW}(Q)$.~$\Box$

\section{The image of the Lusztig cone $L_{st}(Q)$ under $E$ is $L_{PBW}(Q)$}

We show in this section that the image under $E$
of the Lusztig cone $L_{st}(Q)$ is indeed $L_{PBW}(Q)$. It will follow from
results in the next section that
$S_{\mathbf{i}(Q)}^{\mathbf{k}}(L_{st}(Q))=L_{PBW}(Q)$, since we shall see
that $S_{\mathbf{i}(Q)}^{\mathbf{k}}$ and $E$ are
identical on $L_{st}(Q)$.

We suppose that $a=(a_1,a_2,\ldots ,a_N)\in L_{st}(Q)$, and let
$\mathbf{c}=(c_{ij})=E(\mathbf{a})$. We will show that $\mathbf{a}\in
L_{st}(Q)$ if and only if $\mathbf{c}\in L_{PBW}(Q)$.
We translate the linear inequalities defining $L_{st}(Q)$ using the
linear function $E$, and show that they become the defining inequalities of
$L_{PBW}(Q)$.
We also show that $a_i\geq 0$ for all $i$ if and only if $c_{ij}\geq 0$ for
all $1\leq i\leq j\leq n$.
Since $E$ is a linear function, it will follow that $E(L_{st}(Q))=
L_{PBW}(Q)$.

Recall that:
$$\mathbf{i}(Q)=(l_1\searrow 1)(l_2\searrow 1)\cdots (l_a\searrow 1)
(n\searrow 1)
(n\searrow n+1-r_b)(n\searrow n+1-r_{b-1})\cdots (n\searrow n+1-r_1)$$
is compatible with $Q$ (see the end of Section~\ref{compatible}).

The inequalities defining $L_{st}(Q)$ arise from pairs of equal simple
reflections occurring in $\mathbf{i}(Q)$. Note that, for the first $a+1$
factors appearing in the above, each factor is always contained in the one
immediately to the right, and that, for the last $b+1$ factors, each factor
is always contained in the one immediately to the left. It is clear that the
defining inequalities for $L_{st}(Q)$ always arise from such pairs of
factors. To make the notation clearer, let us define $l_{a+1}=r_{b+1}=n$.

Let us consider such a pair, $(l_p\searrow 1)(l_{p+1}\searrow 1)$, and
suppose that $1\leq s\leq l_p$, so $s$ occurs both in $(l_p\searrow 1)$ and
in $(l_{p+1}\searrow 1)$. Suppose $s$ appears in position $t$ of $\mathbf{i}$
in the first factor; it then must appear in the second factor --- suppose that
this is in position $u$ in the second factor. Let us first assume
that $s>1$. The defining inequality of $L_{st}(Q)$ arising from this pair of
$s$'s is:
\begin{equation}
a_{t+1}+a_{u-1}\geq a_t+a_u.
\end{equation}
We can rewrite this as:
\begin{equation}\label{lc}
a_{t}-a_{t+1}\leq a_{u-1}-a_u.
\end{equation}
Since $D(\mathbf{c})=\mathbf{a}$, we have
(from equation~(\ref{Ddef}) in Section~\ref{representations}):
\begin{eqnarray*}
a_t & = & \sum_{\alpha\in T_{n+1-p},\ s\in \alpha}c_{\alpha}
\\
a_{t+1} & = & \sum_{\alpha\in T_{n+1-p}),\ s-1\in
\alpha}c_{\alpha} \\
a_{u-1} & = & \sum_{\alpha\in T_{n-p},\ s+1\in
\alpha}c_{\alpha} \\
a_{u} & = & \sum_{\alpha\in T_{n-p},\ s\in
\alpha}c_{\alpha}
\end{eqnarray*}
Suppose first, that $s\leq p$. Then $s-1\in\alpha$ implies that $s\in
\alpha$ for $\alpha\in T_{n+1-p}$, and $s\in \alpha$ implies that
$s+1\in \alpha$ for $\alpha\in T_{n-p}$, as $s+1\leq p+1$ also.
Thus the inequality (\ref{lc}) becomes:
$$\sum_{\alpha\in T_{n+1-p},s\in
\alpha,s-1\not\in\alpha}c_{\alpha} \leq
\sum_{\alpha\in T_{n-p},s+1\in
\alpha,s\not\in\alpha}c_{\alpha}.$$
This is an inequality of type (L2) or (L3) for $L_{PBW}(Q)$,
depending on what $s$ is.

Next suppose that $s=1$. Then the defining inequality (\ref{lc}) is
the same, except that the term $a_{t+1}$ does not appear. So we have
$a_t\leq a_{u-1}-a_u$. This translates to
$$\sum_{\alpha\in T_{n+1-p},1\in \alpha}c_{\alpha}\leq
\sum_{\alpha\in T_{n-p},2\in
\alpha,1\not\in\alpha}c_{\alpha},$$
and we again get an inequality of type (L2) or (L3) for $L_{PBW}(Q)$.

Suppose next, that $s>p$. Then $s\in\alpha$ implies that $s-1\in
\alpha$ for $\alpha\in T_{n+1-p}$, and $s+1\in \alpha$ implies that
$s\in \alpha$ for $\alpha\in T_{n-p}$, as $s+1\leq p+1$ also.
We rewrite inequality (\ref{lc}) as $a_{t+1}-a_t\geq a_u-a_{u-1}$, and we
have:
$$\sum_{\alpha\in T_{n+1-p},s-1\in
\alpha,s\not\in\alpha}c_{\alpha} \geq
\sum_{\alpha\in T_{n-p},s\in
\alpha,s+1\not\in\alpha}c_{\alpha}.$$
It is clear that this is an inequality of type (L1) or (L4) for $L_{PBW}(Q)$,
depending on what $s$ is. We do not get a boundary case to consider in this
case.

For consecutive factors appearing after $(n\searrow 1)$, a similar argument
shows that they give defining inequalities for $L_{PBW}(Q)$. It is easy to see that,
if all possible pairs of consecutive factors are taken, we get precisely the
defining inequalities for $L_{PBW}(Q)$.

Thus, the defining inequalities for $L_{st}(Q)$ correspond, under $E$,
to the defining inequalities for $L_{PBW}(Q)$. It remains to check that
$a_i\geq 0$ for all $i$ if and only if $c_{ij}\geq 0$ for
all $1\leq i\leq j\leq n$ (if $E(\mathbf{a})=\mathbf{c}$).

By the description of the coordinates of the function $D=E^{-1}$
(see equation~(\ref{Ddef}) in Section~\ref{representations}) as
nonnegative combinations
of the coordinate entries of $\mathbf{c}$, it is clear that if all
$c_{ij}\geq 0$ then all $a_i\geq 0$. Now, suppose that all $a_i\geq 0$.
Then we know by~\cite[4.1]{marsh8}
that, if $\mathbf{a}$ is labelled $(a_{\alpha})_{\alpha\in
\Phi^+}$ according to the ordering on the positive roots induced by
$\mathbf{i}(Q)$, then, whenever $\alpha,\beta\in \Phi^+$ and
$\alpha\geq \beta$ (i.e. $\alpha=\beta$ plus a nonnegative combination
of simple roots), we have $a_{\alpha}\geq a_{\beta}$. By the description of
the function $E$, we know that the
entries on $c_{ij}$, regarded as functions of the coordinates of $\mathbf{a}$,
are all of the form $a_{\alpha}-a_{\beta}$ where $\alpha=\beta+\gamma$, for
some simple root $\gamma$. It follows that $c_{ij}\geq 0$ for all
$1\leq i\leq j\leq n$.

We have proved:

\begin{theorem} \label{conedescription}
$E(L_{st}(Q))=L_{PBW}(Q)$.
\end{theorem}

{\bf Example:} We return to our running example, with $Q=RLRL$ in type
$A_5$. Let $\mathbf{a}\in \mathbb{N}^N$. Then $E(\mathbf{a})=(\mathbf{c})=
(c_{ij})$ where $\mathbf{c}$ is given by the diagram in
Figure~\ref{triangle2}.

\begin{figure}[htbp]

\beginpicture

\setcoordinatesystem units <1cm,1cm>             
\setplotarea x from -8 to 5, y from -0.3 to 4.3       

\put {$a_2-a_1$} at -4 4       %
\put {$a_5-a_4$} at -2 4      %
\put {$a_9-a_8$} at 0 4      %
\put {$a_{13}-a_{12}$} at 2 4      %
\put {$a_{15}$} at 4 4      %
\put {$a_1$} at -3 3      %
\put {$a_4-a_6$} at -1 3      %
\put {$a_8-a_{10}$} at 1 3      %
\put {$a_{12}-a_{14}$} at 3 3       %
\put {$a_6-a_3$} at -2 2      %
\put {$a_{10}-a_7$} at 0 2      %
\put {$a_{14}$} at 2 2      %
\put {$a_3$} at -1 1      %
\put {$a_7-a_{11}$} at 1 1      %
\put {$a_{11}$} at 0 0      %
\endpicture
\caption{The function $S_{\mathbf{i}(Q)}^{\mathbf{k}}(\mathbf{a})$}
\label{triangle2}
\end{figure}

We give below the correspondence between
the defining inequalities of $L_{st}(Q)$ and those of $L_{PBW}(Q)$.
Recall that $\mathbf{i}(Q)=(2,1,4,3,2,1,5,4,3,2,1,5,4,3,5)$.

\begin{tabular}{c|c}
Inequality of $L_{st}(Q)$ & Inequality of $L_{PBW}(Q)$ \\
\hline
$a_2+a_4\geq a_1+a_5$               & $c_{11}\geq c_{22}$ \\
$a_5\geq a_2+a_6$                   & $c_{22}+c_{23}\geq c_{11}+c_{12}$ \\
$a_4+a_7\geq a_3+a_8$               & $c_{23}+c_{13}\geq c_{34}+c_{24}$ \\
$a_5+a_8\geq a_4+a_9$               & $c_{22}\geq c_{33}$ \\
$a_6+a_9\geq a_5+a_{10}$            & $c_{33}+c_{34}\geq c_{22}+c_{23}$ \\
$a_{10}\geq a_6+a_{11}$             & $c_{24}+c_{25}\geq c_{13}+c_{14}$ \\
$a_8\geq a_7+a_{12}$                & $c_{34}+c_{24}\geq c_{45}+c_{35}$ \\
$a_9+a_{12}\geq a_8+a_{13}$         & $c_{33}\geq c_{44}$ \\
$a_{10}+a_{13}\geq a_9+a_{14}$      & $c_{44}+c_{45}\geq c_{33}+c_{34}$ \\
$a_{13}\geq a_{12}+a_{15}$          & $c_{44}\geq c_{55}$
\end{tabular}

\section{Description of the function $S_{\mathbf{i}(Q)}^{\mathbf{k}}$}
\label{description}

We now show that the function $S=S_{\mathbf{i}(Q)}^{\mathbf{k}}$ coincides
with the function $E=D^{-1}$ arising from
representations of quivers, as described in Section~\ref{representations},
on the degeneration cone $C_{st}(Q)$.
We first of all note that each edge of $Q$ corresponds in a natural way to
a slice of $A(Q_{\mathbf{k}})$; such a slice generates the factor of
$\mathbf{i}(Q)$ corresponding to this edge --- see Section~\ref{compatible}.

We suppose that $\mathbf{a}\in C_{st}(Q)=D(C_{PBW}(Q))$
(the degeneration cone).
We will show that $\mathbf{a}\in X_{st}(Q)$ (the string cone), and that
$S(\mathbf{a})=E(\mathbf{a})$. It then follows that, as $E$ and $S$ are
bijective, $S^{-1}(\mathbf{c})=\mathbf{a}=D(\mathbf{c})$ for any
$\mathbf{c}\in C_{PBW}(Q)$. This is the result
we would like to prove, and this section is mainly devoted to achieving this
aim. In the following analysis, $\mathbf{a}$ shall denote an arbitrary
element of $C_{st}(Q)$, and $\mathbf{c}=E(\mathbf{a})$.

In applying $S_{\mathbf{i}(Q)}^{\mathbf{k}}$ to $\mathbf{a}=(a_{ij})$, we need
to compute $\mathbf{c}$ such that $\TF_{i_1}^{a_1}\cdots \TF_{i_k}^{a_k}
\cdot 1 \equiv F_{\mathbf{k}}^{\mathbf{c}}\mod v\mathcal{L}$, where
$\mathbf{i}=\mathbf{i}(Q)$. Recall (see the end of Section~\ref{compatible})
that
$$\mathbf{i}(Q)=(l_1\searrow 1)(l_2\searrow 1)\cdots (l_a\searrow 1)
(n\searrow 1)
(n\searrow n+1-r_b)(n\searrow n+1-r_{b-1})\cdots (n\searrow n+1-r_1).$$

Let $m$ be an edge of $Q$. Then in the Berenstein-Fomin-Zelevinsky
arrangement corresponding to $Q$
(see, for example, Figure~\ref{bfzdiagram}),
the line corresponding to $m$ passes first
through the lines corresponding to $R$'s to the left of $m$, starting from
the right, then the line from top left to bottom right, followed by the $L$'s
to the left of $E$ in $Q$, from left to right. This is how the reduced
expression $\mathbf{i}(Q)$ is built up. Thus we apply the product
\begin{equation}\label{product}
\TF(l_1)\TF(l_2)\cdots \TF(l_a)\TF(n)\TF(r_b)\TF(r_{b-1})\cdots \TF(r_1)
\end{equation}
to $1=F_{\mathbf{k}}^{\mathbf{0}}$ where $\mathbf{0}$ denotes the zero vector,
and $\TF(m)$ is a monomial of Kashiwara operators defined as follows.
Given $m\in [1,n]$, let $d$ be maximal so that $l_d<m$ and let $e$ be maximal
so that $r_e<m$. Then we have that $d+e=m-1$, so $m-e-d=1$, and
$n-e-d=n-m+m-e-d=n+1-m$. Then, if edge $m$ is an $L$, we have:
$$\TF(m)=
\TF_{m}^{a_{r_e+1,m+1}}
\TF_{m-1}^{a_{r_{e-1}+1,m+1}}  \cdots 
\TF_{m-e+1}^{a_{r_1+1,m+1}}
\TF_{m-e}^{a_{1,m+1}}
\TF_{m-e-1}^{a_{l_1+1,m+1}}
\TF_{m-e-2}^{a_{l_2+1,m+1}}  \cdots
\TF_{m-e-d}^{a_{l_d+1,m+1}},
$$
and if edge $m$ is an $R$, or $m=n$, we have:
$$\TF(m)=
\TF_{n}^{a_{r_e+1,m+1}}
\TF_{n-1}^{a_{r_{e-1}+1,m+1}}  \cdots
\TF_{n-e+1}^{a_{r_1+1,m+1}}
\TF_{n-e}^{a_{1,m+1}}
\TF_{n-e-1}^{a_{l_1+1,m+1}}
\TF_{n-e-2}^{a_{l_2+1,m+1}}  \cdots
\TF_{n-e-d}^{a_{l_d+1,m+1}}.
$$

We can regard the vector $\mathbf{c}\in \mathbb{N}^N$ as the vector
$(c_{\alpha})_{\alpha\in\Phi^+}$ (see section~\ref{monomialbases}).
Each root $\alpha\in\Phi_+$ will lie in a slice $T_z$ of $\Phi_+$.
Given $z\in\mathbb{N}$, let $\mathbf{c}(z)$ denote the vector
$\mathbf{c}$ but with $c_{\alpha}$ set to zero for every $\alpha$ lying
in a slice $T_{z'}$ with $z'<z$; thus $\mathbf{c}(n+1)=\mathbf{0}$ and
$\mathbf{c}(1)=\mathbf{c}$.

For $i=1,2,\ldots ,n$, let $\TP_i$ denote the $i$th product appearing
in~(\ref{product}). Thus, for $i=1,2,\ldots ,a$, $\TP_i=\TF(l_i)$,
$\TP_{a+1}=\TF(n)$, and for $i=a+2,\ldots ,n$, $\TP_i=\TF(r_{a+2-i+b})$,
and we have
$$\TP_1\TP_2\cdots \TP_n=
\TF(l_1)\TF(l_2)\cdots \TF(l_a)\TF(n)\TF(r_b)\TF(r_{b-1})\cdots \TF(r_1).$$
We will show that, for $z=1,2,\ldots ,n$,
$\TP_z \cdot F_{\mathbf{k}}^{\mathbf{c(z+1)}}\equiv
F_{\mathbf{k}}^{\mathbf{c(z)}}\mod v\mathcal{L}$, from which it will follow
that $$\TP_1\TP_2\cdots \TP_n\cdot 1\equiv F_{\mathbf{k}}^{\mathbf{c}}$$
as required.

We also need to show that $\mathbf{a}\in X_{st}(Q)$. We will use the
following definiton: A monomial action
$\TF_{j_1}^{b_1}\TF_{j_2}^{b_2}\cdots \TF_{j_t}^{b_t}\cdot
F_{\mathbf{k}}^{\mathbf{x}}$
is said to satisfy
(STRING) provided that
$$\TE_{j_u}\TF_{j_{u+1}}^{b_{u+1}}\TF_{j_{u+2}}^{b_{u+2}}\cdots \TF_{j_{t}}^{b_t}
\cdot F_{\mathbf{k}}^{\mathbf{x}}\equiv 0\mod v\mathcal{L},$$ for
$u=1,2,\ldots ,t$.
We will also use the description of the action of the $\TF_i$'s on a PBW basis
modulo $v\mathcal{L}$ as proved by the second author in~\cite{reineke1}:

\begin{prop} \label{reineke} (Reineke)
Suppose ${\bf c}=(c_{ij})\in \mathbb{N}^N$.
For each $1\leq i\leq j\leq n$, define
$$f_{ij}=\sum_{k=1}^{i} c_{kj}-\sum_{k=1}^{i-1} c_{k,j-1}.$$
Let $i_0$ be maximal so that $f_{i_0j}=\max_i f_{ij}$. Then $\widetilde{F}_j$
increases $c_{i_0j}$ by $1$, decreases $c_{i_0,j-1}$ by $1$ (unless
$i_0=j$, when this latter effect does not occur), and leaves the other
$c_{ij}$'s unchanged.
\end{prop}

It is easy to see, using Theorem~\ref{kashiwara}(ii), that this implies the
following description of the $\TE_i$'s:

\begin{prop} \label{reineke2}
Suppose ${\bf c}=(c_{ij})\in \mathbb{N}^N$.
For each $1\leq i\leq j\leq n$, define
$$f_{ij}=\sum_{k=1}^{i} c_{kj}-\sum_{k=1}^{i-1} c_{k,j-1}.$$
Let $i_0$ be minimal so that $f_{i_0j}=\max_i f_{ij}$. Then $\widetilde{E}_j$
decreases $c_{i_0j}$ by $1$, inreases $c_{i_0,j-1}$ by $1$ (unless
$i_0=j$, when this latter effect does not occur), and leaves the other
$c_{ij}$'s unchanged, except, if $c_{i_0j}=0$, then it acts as zero modulo
$q\mathcal{L}$.
\end{prop}

We will need the following technical Lemmas, describing the action of
$\TF_i$ and $\TE_i$ in certain circumstances:

\begin{lemma}\label{technical1}
Fix $1\leq i\leq j\leq n$ and $s\in\N$. Suppose the following hold for a 
triangle $(c_{kl})\in \mathbb{N}^N$:

(a) For $k=1\ldots i-1$, we have 
$\sum_{l=k+1}^ic_{lj}\geq\sum_{l=k}^{i-1}c_{l,j-1}$. \\
(b) For $k=i+1\ldots j$, we have 
$s\leq\sum_{l=i}^{k-1}c_{l,j-1}-\sum_{l=i+1}^kc_{lj}$.

Then $\TF_j^s$ acts on $(c_{kl})$ by increasing $c_{ij}$ by $s$, decreasing 
$c_{i,j-1}$ by $s$, and leaving the other $c_{kl}$'s unchanged.
\end{lemma}

\noindent {\bf Proof:} For $t=0\ldots s-1$, define a new triangle $c^t_{kl}$ by
$$c^t_{kl}=\left\{\begin{array}{ll}
c_{ij}+t, & k=i,l=j \\
c_{i,j-1}-t, & k=i, l=j-1 \\
c_{kl}, & \mbox{otherwise}. \end{array}\right.$$
The Lemma clearly holds if for all $t=0\ldots i-1$, the index $i_0$ of 
Proposition \ref{reineke} equals $i$. This translates into the following
inequalities for the $f_{kl}$'s of Proposition \ref{reineke}:
$$f_{kj}\leq f_{ij}\mbox{ for $k=1\ldots i-1$},\;\;\; f_{kj}<f_{ij}\mbox{ for 
$k=i+1\ldots j$}.$$
Using the definitions of $(c^t_{kl})$ and $f_{kl}$, this translates into the 
following conditions:
$$t\geq\sum_{l=k}^{i-1}c_{l,j-1}-\sum_{l=k+1}^ic_{lj}\mbox{ for $t=0\ldots 
s-1$, $k=1\ldots i-1$},$$
which is equivalent to condition (a), and
$$t<\sum_{l=i}^{k-1}c_{l,j-1}-\sum_{l=i+1}^kc_{lj}\mbox{ for $t=0\ldots s-1$, 
$k=i+1\ldots j$},$$
which is equivalent to condition (b).~$\Box$

\begin{lemma}\label{technical2}
Fix $1\leq i\leq j\leq n$ and $s\in \mathbb{N}$. Suppose the following hold
for a triangle $(c_{kl})\in \mathbb{N}^N$:

(a) For $k=1\ldots i$, we have $c_{kj}=0$. \\
(b) For $k=i+1\ldots j$, we have $s\leq 
\sum_{l=i}^{k-1}c_{l,j-1}-\sum_{l=i+1}^kc_{lj}$.

Then for $\widetilde{s}:=s+\sum_{k=1}^{i-1}c_{k,j-1}$, we have:

$${\TF_j^{\widetilde{s}}(c_{kl})}_{kl}=\left\{ \begin{array}{ll}
0, & l=j-1,k<i \\
c_{k,j-1}, & l=j,k<i \\
c_{i,j-1}-s, & l=j-1,k=i \\
s, & l=j, k=i \\
c_{kl}, & \mbox{otherwise}.
\end{array} \right.
$$
\end{lemma}

\noindent {\bf Proof:} We proceed by induction on $i$. If $i=1$, the claimed 
statement follows directly from Lemma \ref{technical1}. For arbitrary $i$, we 
set $i'=i-1$, $s'=c_{i-1,j-1}$ and claim that we can apply the Lemma -- which 
we assume to be already true for $i'$ -- with the same triangle $(c_{kl})$, 
but with $(i,j,s)$ replaced by $(i',j,s')$. We thus have to show that the 
assumptions of the Lemma are satisfied:\\
Condition (a) is satisfied trivially. For $k=i+1\ldots j$, condition (b) gives 
$\sum_{l=i}^{k-1}c_{l,j-1}-\sum_{l=i+1}^kc_{lj}\geq s\geq 0$, which (using 
$c_{ij}=0$ by condition (a)) means
$$\sum_{l=i'}^{k-1}c_{l,j-1}-\sum_{l=i'+1}^kc_{lj}\geq c_{i-1,j-1}=s'.$$
For $k=i'+1=i$, the desired condition is $s'\leq c_{i-1,j-1}-c_{ij}$, which 
clearly holds since $c_{ij}=0$.\\
Thus, we can apply the Lemma to $(i',j,s')$. Setting 
$(c'_{kl})=\TF_j^{\widetilde{s'}}(c_{kl})$, we find:
$$c'_{kl}=\left\{\begin{array}{ll}
0,&l=j-1,k<i'\\ c_{k,j-1},&l=j,k<i'\\ c_{i',j-1}-s,&l=j-1,k=i'\\ s,&l=j,k=i'\\ 
c_{kl},&\mbox{otherwise}\end{array}\right\}=\left\{\begin{array}{ll}0,&l=j-1,k<i\\
c_{k,j-1},&l=j,k<i\\
c_{kl},&\mbox{otherwise}.\end{array}\right.$$
Since $\widetilde{s'}=\sum_{k=1}^{i-1}c_{k,j-1}$, it thus remains to show that 
$\TF_j^s$ acts on $(c'_{kl})$ by decreasing $c_{i,j-1}$ by $s$, increasing 
$c_{ij}=0$ by $s$, and leaving the rest of $(c'_{kl})$ unchanged. To prove 
this, we only have to check the assumptions of Lemma \ref{technical1} for 
$(c'_{kl})$: for condition (a) of Lemma \ref{technical1}, this is trivial 
since $c'_{k,j-1}=0$ for $k<i$; for condition (b) of Lemma \ref{technical1}, 
we just use condition (b) of the present Lemma. Applying Lemma 
\ref{technical1}, we see that we are done.~$\square$

\begin{lemma} \label{technical3}
Fix $1\leq i\leq j\leq n$ and $s\in \mathbb{N}$. Suppose the following hold
for a triangle $(c_{kl})\in \mathbb{N}^N$:

(a) For $k=1\ldots i$, we have $c_{kj}=0$. \\
(b) For $k=i+1\ldots j$, we have $0\leq 
\sum_{l=i}^{k-1}c_{l,j-1}-\sum_{l=i+1}^kc_{lj}$.

Then we have $\TE_j(c_{kl})=0$.
\end{lemma}

{\bf Proof:} For each $1\leq p\leq j\leq n$, define
$$f_{pj}=\sum_{k=1}^{p} c_{kj}-\sum_{k=1}^{p-1} c_{k,j-1},$$
as in Proposition~\ref{reineke2}.
By assumption (a), $f_{pj}=-\sum_{k=1}^{p-1} c_{k,j-1}$, for $p=1,2,\ldots i$,
and, using both assumptions,
\begin{eqnarray*}
f_{pj} & = & \sum_{k=i+1}^p c_{kj} -\sum_{k=1}^{p-1} c_{k,j-1} \\
       & = & f_{ij}+\sum_{k=i+1}^p c_{kj} -\sum_{k=i}^{p-1} c_{k,j-1}\leq f_{ij}
\end{eqnarray*}
for $p=i+1,i+2,\ldots ,j$.
Hence if $p_0$ is minimal so that $f_{p_0j}=\max_p f_{pj}$, we must have
$p_0\leq i$. It follows from (a) that $c_{p_0j}=0$, so we conclude that
$\TE_j$ acts as zero modulo $q\mathcal{L}$.~$\square$

We start by considering the action of the monomial $\TF(r_p)$ (where
$1\leq p\leq b+1$ --- recall that $r_{b+1}=n$) on $\mathbf{c}(n+2-p)$.
Recall that, for $p=1,2,\ldots ,b+1$, we have
$$\TF(r_p)=
\TF_{n}^{a_{r_{p-1}+1,m+1}}
\TF_{n-1}^{a_{r_{p-2}+1,m+1}}  \cdots
\TF_{n-p+2}^{a_{r_1+1,m+1}}
\TF_{n-p+1}^{a_{1,m+1}}
\TF_{n-p}^{a_{l_1+1,m+1}}
\TF_{n-p-1}^{a_{l_2+1,m+1}}  \cdots
\TF_{n+1-p-d}^{a_{l_d+1,m+1}}.
$$

The initial part of the computation is reasonably easy:

\begin{lemma} \label{firstfactor}
Suppose that $1\leq p\leq b$.
Let $\mathbf{d}\in \mathbb{N}^N$ be the triangle given by
\begin{eqnarray*}
d_{n+1-p,n+1-p} & = & a_{1,r_p+1}-a_{l_1+1,r_p+1}, \\
d_{n-p,n+1-p} & = & a_{l_1+1,r_p+1}-a_{l_2+1,r_p+1}, \\
d_{n-p-1,n+1-p} & = & a_{l_2+1,r_p+1}-a_{l_3+1,r_p+1}, \\
\vdots & = & \vdots \\
d_{n+2-p-d,n+1-p} & = & a_{l_{d-1}+1,r_p+1}-a_{l_d+1,r_p+1}, \\
d_{n+3-p-d,n+1-p} & = & a_{l_d+1,r_p+1}.
\end{eqnarray*}
Then we have
$$\TF_{n-p+1}^{a_{1,m+1}}
\TF_{n-p}^{a_{l_1+1,m+1}}
\TF_{n-p-1}^{a_{l_2+1,m+1}}  \cdots
\TF_{n+1-p-d}^{a_{l_d+1,m+1}} \cdot
F_{\mathbf{k}}^{\mathbf{c}(n+2-p)}\equiv
F_{\mathbf{k}}^{\mathbf{c}(n+2-p)+\mathbf{d}}.$$
Furthermore, the action
$$\TF_{n-p+1}^{a_{1,m+1}}
\TF_{n-p}^{a_{l_1+1,m+1}}
\TF_{n-p-1}^{a_{l_2+1,m+1}}  \cdots
\TF_{n+1-p-d}^{a_{l_d+1,m+1}} \cdot
F_{\mathbf{k}}^{\mathbf{c}(n+2-p)}$$
satisfies (STRING).
\end{lemma}

{\bf Proof:}
We first note that we have $a_{1,r+1}\geq a_{l_1+1,r+1}\geq a_{l_2+1,r+1}\geq
\cdots \geq a_{l_b+1,r+1}$. These follow from the fact that
$\mathbf{a}=D(\mathbf{c})$, since we know that all the coordinates of
$\mathbf{c}$ are nonnegative: we use equation~(\ref{Ddef}) in
Section~\ref{representations}.
Let $l_0=0$. We know that, for $t=0,1,\ldots d-1$,
\begin{equation}\label{sum1}
a_{l_t,r+1}=\sum_{\alpha\in T_{n+1-p},n-p+1-t\in\alpha}c_{\alpha},
\end{equation}
and
\begin{equation}\label{sum2}
a_{l_{t+1},r+1}=\sum_{\alpha\in T_{n+1-p},n-p-t\in\alpha}c_{\alpha}.
\end{equation}
Let $t\in \{0,1,\ldots ,d-1\}$.
Since $n-p+1-t\leq n-p+1$, it follows that if $\alpha_{n-p-t}$ appears
in a root of slice $p$, so does $\alpha_{n-p+1-t}$. Thus, by equations
(\ref{sum1}) and (\ref{sum2}), we have $a_{l_t,r+1}\geq a_{l_{t+1},r+1}$
as required.

Let $\mathbf{x}=\mathbf{c}(n+2-p)$.
It is easy to see, using Proposition~\ref{reineke}, that
$\TF_{n+1-r_p}^{a_{l_b+1,r_p+1}}$ sets $x_{n+1-r_p,n+1-r_p}$ to be
$a_{l_b+1,r_p+1}$,
and leaves the other $x_{ij}$'s unchanged. Similarly, it is easy to
see, using Proposition~\ref{reineke2}, that $\TE_{n+1-r_p}$ acts as zero on
$\mathbf{x}$. In both cases this is because $x_{k,n+1-r_p}=0$ for
$k=1,2,\ldots n+1-r$, and $x_{k,n-r_p}=0$ for $k=1,2,\ldots n-r_p$, so
that all of the $f_{ij}$ in Proposition~\ref{reineke} or
Proposition~\ref{reineke2} are zero.

Next,
$\TF_{n+2-r_p}^{a_{l_{b-1}+1,r_p+1}}$ sets $x_{n+1-r_p,n+1-r_p}=0$
and sets $x_{n+1-r_p,n+2-r_p}=a_{l_b+1,r_p+1}$.
It also sets $x_{n+2-r_p,n+2-r_p}$ to be $a_{l_{b-1}+1}-a_{l_b+1,r_p+1}$. The
rest of the $\TF_i$'s act in a similar way, until we finally get the
vector described in the Lemma. It is also easy to see that the given
monomial action satisfies (STRING) using Proposition~\ref{reineke2}.~$\Box$.

The next part of the computation, which involves computing
$\TF_{n}^{a_{r_{p-1}+1,m+1}}
\TF_{n-1}^{a_{r_{p-2}+1,m+1}}  \cdots
\TF_{n-p+2}^{a_{r_1+1,m+1}} \cdot
F_{\mathbf{k}}^{\mathbf{c'}(n-p+2)}$
(where $\mathbf{c'}(n-p+2)=\mathbf{c}(n-p+2)+\mathbf{d}$ where $\mathbf{d}$ is
as in Lemma~\ref{firstfactor}), is more involved. We need to consider what
happens step-by-step, and to understand what is happening in detail along
each slice $T_z=T_{n+1-p}$ of $A(Q_{\mathbf{k}})$.
Let $T_z=\{\beta_{z1},\beta_{z2},\ldots \beta_{zk_z}\}$, where the $k_z$
roots in $T_z$ are listed according to increasing height. Thus $\beta_{z1}=
\alpha_{n+1-p}$.

Given $q$ such that $1\leq q\leq p-1$, we define a vector $\mathbf{c}(p,q)$
as follows. (In order to simplify notation, we use the numbering arising
from the edges of $Q$ oriented to the right, rather than the usual slice
numbering.)
If $\alpha$ belongs to a slice numbered $n+2-p$ or greater (i.e.
corresponding to one of the right-oriented edges numbered $r_1,r_2,\ldots ,
r_{p-1}$), then
$c_{\alpha}(p,q)=c_{\alpha}$. We also set, for
$i=1,2,\ldots ,r_{q}-1$, $c_{\beta_{zi}}(p,q)=a_{i,r_p+1}-a_{i+1,r_p+1}$,
and set $c_{\beta_{r_q,i}}(p,q)=a_{r_q+1,r_p+1}-a_{l_s+1,r_p+1}$, where
$l_t$ is the number of the edge of the first $L$ to the right of edge $r_q$.
Suppose that $\beta_{zr_{q}}=\alpha_{ij}$. Then we also set
$c_{i-1,j}(p,q)=a_{l_t+1,r_p+1}-a_{l_{t+1}+1,r_p+1}$,
$c_{i-2,j}(p,q)=a_{l_{t+2}+1, r_p+1},\ldots ,c_{t+i+1-b,j}(p,q)=
a_{l_{b-1}+1,r_p+1}-a_{l_b+1,r_p+1}$, and
$c_{t+i-b,j}(p,q)=a_{l_b+1,r_p+1}$. We write $\mathbf{c}(p,0)$ for
the vector $\mathbf{c}(n+2-p)+\mathbf{d}$ appearing in
Lemma~\ref{firstfactor}. The vector $\mathbf{c}(p,q-1)$ (which
appears in the Lemma below) can be visualised as in
Figure~\ref{slicediagram} (note that the case displayed is where
edge $r_q-1$ is oriented to the right).
\begin{figure}[htbp]
\beginpicture

\setcoordinatesystem units <0.8cm,0.8cm>             
\setplotarea x from -7 to 5, y from -1 to 21       
\linethickness=0.5pt           


\put{$a_{l_b+1,r+1}$} at 0 1 %
\put{$a_{l_{t+1}+1,r+1}-a_{l_{t+2}+1,r+1}$} at 2 3 %
\put{$a_{l_{t}+1,r+1}-a_{l_{t+1}+1,r+1}$} at 3 4 %
\put{$a_{r_{q}-1+1,r+1}-a_{l_{t}+1,r+1}$} at 4 5 %
\put{$a_{l_{t-1}+2,r+1}-a_{l_{t-1}+3,r+1}$} at 2 7 %
\put{$a_{l_{t-1}+1,r+1}-a_{l_{t-1}+2,r+1}$} at 1 8 %
\put{$a_{l_{t-1},r+1}-a_{l_{t-1}+1,r+1}$} at 2 9 %
\put{$a_{l_{3},r+1}-a_{l_{3}+1,r+1}$} at 2 11 %
\put{$a_{l_{2}+1,r+1}-a_{l_{2}+2,r+1}$} at 0 13 %
\put{$a_{l_{2},r+1}-a_{l_{2}+1,r+1}$} at 1 14 %
\put{$a_{l_{1}+2,r+1}-a_{l_{1}+3,r+1}$} at -1 16 %
\put{$a_{l_{1}+1,r+1}-a_{l_{1}+2,r+1}$} at -2 17 %
\put{$a_{l_{1},r+1}-a_{l_{1}+1,r+1}$} at -1 18 %
\put{$a_{2,r+1}-a_{3,r+1}$} at -3 20 %
\put{$a_{1,r+1}-a_{2,r+1}$} at -4 21 %


\setlinear \plot 0.2 1.2  0.6 1.6 / %
\setdashes <0.8mm,0.8mm>          %
\setlinear \plot 0.6 1.6  1.4 2.4 / %
\setsolid
\setlinear \plot 1.4 2.4  1.8 2.8 / %

\setlinear \plot 2.2 3.2  2.8 3.8 / %

\setlinear \plot 3.2 4.2  3.8 4.8 / %

\setlinear \plot 3.8 5.2  3.4 5.6 / %
\setdashes <0.8mm,0.8mm>          %
\setlinear \plot 3.4 5.6  2.6 6.4 / %
\setsolid
\setlinear \plot 2.6 6.4  2.2 6.8 / %

\setlinear \plot 1.8 7.2  1.2 7.8 / %

\setlinear \plot 1.2 8.2  1.8 8.8 / %

\setdashes <1mm, 1mm>
\setlinear \plot 2 9.2    2 10.8 / %
\setsolid

\setlinear \plot 1.8 11.2  1.4 11.6 / %
\setdashes <0.8mm,0.8mm>          %
\setlinear \plot 1.4 11.6  0.6 12.4 / %
\setsolid
\setlinear \plot 0.6 12.4  0.2 12.8 / %

\setlinear \plot 0.2 13.2  0.8 13.8 / %

\setlinear \plot 0.8 14.2  0.4 14.6 / %
\setdashes <0.8mm,0.8mm>          %
\setlinear \plot 0.4 14.6  -0.4 15.4 / %
\setsolid
\setlinear \plot -0.4 15.4  -0.8 15.8 / %

\setlinear \plot -1.2 16.2  -1.8 16.8 / %

\setlinear \plot -1.8 17.2  -1.2 17.8 / %

\setlinear \plot -1.2 18.2  -1.6 18.6 / %
\setdashes <0.8mm,0.8mm>          %
\setlinear \plot -1.6 18.6  -2.4 19.4 / %
\setsolid
\setlinear \plot -2.4 19.4  -2.8 19.8 / %

\setlinear \plot -3.2 20.2  -3.8 20.8 / %


\put{\bf\large A} at 1.5 5


\setcoordinatesystem point at -4 0

\put{\bf\large B} at 3.6 4

\put{$0$} at -1 0 %
\put{$0$} at 1 2 %
\put{$0$} at 2 3 %
\put{$0$} at 3 4 %

\setlinear \plot -0.8 0.2  -0.4 0.6 / %
\setdashes <0.8mm,0.8mm>          %
\setlinear \plot -0.4 0.6  0.4 1.4 / %
\setsolid
\setlinear \plot 0.4 1.4  0.8 1.8 / %

\setlinear \plot 1.2 2.2  1.8 2.8 / %

\setlinear \plot 2.2 3.2  2.8 3.8 / %

\setlinear \plot 3.2 4.2  3.8 4.8 / %

\setlinear \plot 3.8 5.2  3.4 5.6 / %
\setdashes <0.8mm,0.8mm>          %
\setlinear \plot 3.4 5.6  2.6 6.4 / %
\setsolid
\setlinear \plot 2.6 6.4  2.2 6.8 / %

\setlinear \plot 1.8 7.2  1.2 7.8 / %

\setlinear \plot 1.2 8.2  1.8 8.8 / %

\setdashes <1mm, 1mm>
\setlinear \plot 2 9.2    2 10.8 / %
\setsolid

\setlinear \plot 1.8 11.2  1.4 11.6 / %
\setdashes <0.8mm,0.8mm>          %
\setlinear \plot 1.4 11.6  0.6 12.4 / %
\setsolid
\setlinear \plot 0.6 12.4  0.2 12.8 / %

\setlinear \plot 0.2 13.2  0.8 13.8 / %

\setlinear \plot 0.8 14.2  0.4 14.6 / %
\setdashes <0.8mm,0.8mm>          %
\setlinear \plot 0.4 14.6  -0.4 15.4 / %
\setsolid
\setlinear \plot -0.4 15.4  -0.8 15.8 / %

\setlinear \plot -1.2 16.2  -1.8 16.8 / %

\setlinear \plot -1.8 17.2  -1.2 17.8 / %

\setlinear \plot -1.2 18.2  -1.6 18.6 / %
\setdashes <0.8mm,0.8mm>          %
\setlinear \plot -1.6 18.6  -2.4 19.4 / %
\setsolid
\setlinear \plot -2.4 19.4  -2.8 19.8 / %

\setlinear \plot -3.2 20.2  -3.8 20.8 / %

\endpicture
\caption{The vector $\mathbf{c}(p,q-1)$.}
\label{slicediagram}
\end{figure}
We have:
\begin{lemma} \label{secondfactor}
Suppose that $1\leq p\leq b$ and that $1\leq q\leq p-1$. We have that

(a) $\TF_{n+1-p+q}^{a_{r_q+1,r_p+1}}\cdot F_{\mathbf{k}}^{\mathbf{c}(p,q-1)}\equiv
F_{\mathbf{k}}^{\mathbf{c}(p,q)}\mod v\mathcal{L}$, and that \\
(b) $\TE_{n+1-p+q}^{a_{r_q+1,r_p+1}}\cdot F_{\mathbf{k}}^{\mathbf{c}(p,q-1)}\equiv
0\mod v\mathcal{L}$.
\end{lemma}

{\bf Proof:} Suppose that $\beta_{z,r_q-1}=\alpha_{ij}$. It is easy
to check that $i=n+1-p-(t-1)$ and $j=n+1-p+q-1$, where $l_t$ is the
number of the edge of the first $L$ to the right of edge $r_q-1$.
Let $s=a_{r_q+1,r_p+1}-a_{l_t+1,r_p+1}$, 
(Note that these notations were used for $\mathbf{c}(p,q)$ in the above,
rather than $\mathbf{c}(p,q-1)$). 

We apply Lemma~\ref{technical2} with the pair $i,j+1$ and $s$ as above.
We need to check the assumptions of the Lemma.
It is clear from the definition of $\mathbf{c}(p,q-1)$ that
$c_{k,j+1}(p,q-1)=0$ for $k=1,2,\ldots ,i$. Suppose that $i+1\leq k\leq j+1$.
We consider 
$$\sum_{l=i}^{k-1}c_{l,j}(p,q-1)-\sum_{l=i+1}^{k}c_{l,j+1}(p,q-1).$$
This is a sum of terms, each of which is the sum of elements on one
side on an inequality (C1) or (C2) (see section~\ref{mrcones}),
minus the corresponding sum of elements on the other side of the inequality.
Since the vector $\mathbf{c}\in C_{PBW}(Q)$ satisfies
inequalities (C1) and (C2), we have that:
\begin{equation}\label{satisfied1}
\sum_{l=i}^{k-1}c_{l,j}-\sum_{l=i+1}^{k}c_{l,j+1}\geq 0
\end{equation}
(indeed, the inequalities of $C_{PBW}(Q)$ are designed to ensure this holds).
By the definition of $\mathbf{c}(p,q-1)$, it is clear that for
$l=i+1,\ldots j$, $c_{l,j}=c_{l,j}(p,q-1)$ and that for $l=i+1,\ldots ,k$,
$c_{l,j+1}=c_{l,j+1}(p,q-1)$. Furthermore, $c_{ij}(p,q-1)=
a_{r_q-1+1,r_p+1}-a_{l_{t}+1-r_p+1}=
s+(a_{r_q-1+1,r_p+1}-a_{r_q+1,r_p+1})=s+c_{ij}$.
It follows from Equation~(\ref{satisfied1}) that
\begin{equation}\label{satisfied2}
\sum_{l=i}^{k-1}c_{l,j}(p,q-1)-\sum_{l=i+1}^{k}c_{l,j+1}(p,q-1)\geq s
\end{equation}
as required in Lemma~\ref{technical2}. Thus the conditions of the Lemma are
all satisfied, and it applies. Let $\tilde{s}=s+\sum_{k=1}^{i-1}c_{k,j}=
a_{r_{q}+1,r_p+1}-a_{l_t+1,r_p+1}+a_{l_{t}+1,r_p+1}-a_{l_{t+1}+1,r_p+1}+
a_{l_{t+1}+1,r_p+1}-a_{l_{t+2}+1,r_p+1}+\cdots +a_{l_{b-1}+1,r_p+1}-
a_{l_b+1,r_p+1}+a_{l_b+1,r_p+1}=a_{r_q+1,r_p+1}$. Then it is easy to check
that Lemma~\ref{technical2} tells us that
$\TF_{n+1-p+q}^{a_{r_q+1,r_p+1}}\cdot F_{\mathbf{k}}^{\mathbf{c}(p,q-1)}\equiv
F_{\mathbf{k}}^{\mathbf{c}(p,q)}\mod v\mathcal{L}$, giving
us (a) above. This can easily be visualised in Figure~\ref{slicediagram};
the value of $c_{ij}(p,q-1)$ is displayed at $A$ in the diagram, and this
is replaced by $a_{r_{q-1}+1,r_p+1}-a_{r_q+1,r_p+1}$; the value of
$c_{i,j+1}(p,q-1)$ (which is zero) is displayed at $B$ in the diagram, and
this is replaced by $s=a_{r_q+1,r_p+1}-a_{l_t+1,r_p+1}$. Finally, the values
on the diagonal below and to the left of $A$ are all moved one step down
and to the right (to places where the value is currently zero), forming
the diagonal below and to the left of $B$. The new values in the diagonal
below $A$ are all zero.

In order to show (b), we apply Lemma~\ref{technical3} to
$\mathbf{c}(p,q-1)$, again with the pair $i,j+1$. It is clear that
the conditions of this Lemma are satisfied by $\mathbf{c}(p,q-1)$, since
they are same as those in Lemma~\ref{technical2} but with $s$ taken to be
zero. The Lemma tells us that
$\TE_{n+1-p+q}^{a_{r_q+1,r_p+1}}\mathbf{c}(p,q-1)\equiv 0\mod v\mathcal{L}$,
giving us (b) above.~$\Box$

We can now prove our main result:

\begin{theorem} \label{coincide}
Let $\mathbf{c}\in C_{PBW}(Q)$. Then
$(S_{\mathbf{i}(Q)}^{\mathbf{k}})^{-1}(\mathbf{c})=D(\mathbf{c})$.
\end{theorem}

{\bf Proof:}
Let $\mathbf{a}\in C_{st}(Q)=D(C_{PBW}(Q))$, and suppose $1\leq p\leq b$.
Lemma~\ref{firstfactor} tells us that
$$\TF_{n-p+1}^{a_{1,m+1}}
\TF_{n-p}^{a_{l_1+1,m+1}}
\TF_{n-p-1}^{a_{l_2+1,m+1}}  \cdots
\TF_{n+1-p-d}^{a_{l_d+1,m+1}} \cdot
F_{\mathbf{k}}^{\mathbf{c}(n+2-p)}\equiv
F_{\mathbf{k}}^{\mathbf{c}(p,0)}.$$
Lemma~\ref{secondfactor}(a) tells us that, for $1\leq q\leq p-1$, we have
$\TF_{n+1-p+q}^{a_{r_q+1,r_p+1}}\cdot F_{\mathbf{k}}^{\mathbf{c}(p,q-1)}\equiv
F_{\mathbf{k}}^{\mathbf{c}(p,q)}\mod v\mathcal{L}$.
Repeated application of this second result (for $q=1,2,\ldots ,p-1$) tells
us that
$\TF(r_p)\cdot F_{\mathbf{k}}^{\mathbf{c}(n+2-p)}\equiv
F_{\mathbf{k}}^{\mathbf{c}(n+1-p)}\mod v\mathcal{L}$, for $p=1,2,\ldots ,r_b$.
It follows that $\TF(n)\TF(r_b)\TF(r_{b-1})\cdots \TF(r_1)\cdot 
F_{\mathbf{k}}^{\mathbf{c}(n+1)}\equiv
F_{\mathbf{k}}^{\mathbf{c}(a+1)}\mod v\mathcal{L}$.
A proof similar to the one given above
can be used to show that $\TF(l_i)\cdot F_{\mathbf{k}}^{\mathbf{c}(i+1)}\equiv
F_{\mathbf{k}}^{\mathbf{c}(i)}\mod v\mathcal{L}$ for $i=1,2,\ldots ,a$.
It then follows that
$\TF(l_1)\TF(l_2)\cdots \TF(l_a)\TF(n)\TF(r_b)\TF(r_{b-1})\cdots \TF(r_1)
\cdot F_{\mathbf{k}}^{\mathbf{c}(n+1)}\equiv
F_{\mathbf{k}}^{\mathbf{c}(1)}\equiv
F_{\mathbf{k}}^{\mathbf{c}}\equiv
F_{\mathbf{k}}^{E(\mathbf{a})}\mod v\mathcal{L},$
as required. It also follows from Lemmas~\ref{firstfactor}
and~\ref{secondfactor} (and a similar proof for the monomial
$\TF(l_1)\TF(l_2)\cdots \TF(l_a)$) that the monomial action \linebreak
$\TF(l_1)\TF(l_2)\cdots \TF(l_a)\TF(n)\TF(r_b)\TF(r_{b-1})\cdots \TF(r_1)
\cdot F_{\mathbf{k}}^{\mathbf{0}}\equiv F_{\mathbf{k}}^{\mathbf{c}}$ satisfies
(STRING), from which it follows that
$\mathbf{a}\in X_{st}(Q)$.

We thus have that, for all $\mathbf{a}\in C_{st}(Q)=D(C_{PBW}(Q))$,
$S(\mathbf{a})=E(\mathbf{a})$ and $\mathbf{a}\in X_{st}(Q)$. It follows
that for all $\mathbf{c}\in C_{PBW}(Q)$,
$(S_{\mathbf{i}(Q)}^{\mathbf{k}})^{-1}(\mathbf{c})=D(\mathbf{c})$,
as required.~$\Box$

We conjecture that this theorem holds in the case when $\mathbf{k}$ is
replaced by any reduced expression for $w_0$ compatible with a quiver.

We have the corollary:

\begin{cor}
$(S_{\mathbf{i}(Q)}^{\mathbf{k}})^{-1}(C_{PBW}(Q))=C_{st}(Q)$.
\end{cor}

{\bf Proof:} This follows from the Theorem and the definition of the
degeneration cone $C_{st}(Q)=D(C_{PBW}(Q))$.~$\square$

Finally, since (by Theorem~\ref{conedescription}),
$E(L_{st}(Q))=L_{PBW}(Q)$ and since $L_{PBW}(Q)\subseteq C_{PBW}(Q)$, we
have, by Theorem~\ref{coincide}, that
$S$ and $E$ coincide on $L_{st}(Q)$. It follows that:

\begin{theorem} \label{image}
$S_{\mathbf{i}(Q)}^{\mathbf{k}}(L_{st}(Q))=L_{PBW}(Q)$.
\end{theorem}

{\bf Example:} We return to our running example, with $Q=RLRL$ in type
$A_5$. By Theorem~\ref{coincide}, for $\mathbf{a}$ such that
$E(\mathbf{a})\in C_{PBW}(Q)$, we have
$S_{\mathbf{i}(Q)}^{\mathbf{k}}(\mathbf{a})=(\mathbf{c})=(c_{ij})$ where
$c_{ij}$ is given in Figure~\ref{triangle2}.

Using the description of $D$, we have, for $\mathbf{c}\in
L_{PBW}(Q)$, that $(S_{\mathbf{i}(Q)}^{\mathbf{k}})^{-1}(\mathbf{c})=
\mathbf{a}$, where $\mathbf{a}=
(c_{12}, c_{11} + c_{12}, c_{14}, c_{23} + c_{13} + c_{14},
c_{22} + c_{23} + c_{13} + c_{14}, c_{13} + c_{14}, c_{25} + c_{15},
c_{34} + c_{24} + c_{25} + c_{15}, c_{33} + c_{34} + c_{24} + c_{25} + c_{15},
c_{24} + c_{25} + c_{15}, c_{15}, c_{45} + c_{35}, c_{44} + c_{45} + c_{35},
c_{35}, c_{55})$.

{\bf Acknowledgements}

The authors gratefully acknowledge support for this research from a
University of Leicester Research Fund Grant, and would like to thank
the Universities of Wuppertal and Leicester, where this research was carried
out, for their hospitality.

\newcommand{\noopsort}[1]{}\newcommand{\singleletter}[1]{#1}

\end{document}